\documentclass[12pt]{article}
\usepackage{amsmath}
\usepackage{amssymb}
\usepackage[mathscr]{eucal}
\usepackage[usenames]{color}
\usepackage[utf8]{inputenc}
\usepackage{url}
\usepackage{longtable}
\usepackage{graphicx}
\usepackage{epstopdf}
\usepackage{subfigure}

\usepackage{authblk}

\usepackage{url}
\usepackage[
bookmarks=true,
backref=true,
colorlinks=true,
linkcolor=red,
citecolor=blue,
urlcolor=green]{hyperref}

\newtheorem{theorem}{Theorem}
\newtheorem{lemma}{Lemma}

\newtheorem{corollary}{Corollary}
\newtheorem{remark}{Remark}

\numberwithin{equation}{section}

\definecolor{darkolivegreen}{rgb}{0.333333, 0.419608, 0.1843140}

\title{\Large Hopf Bifurcation in a Generalized Goodwin Model with Delay
}

\author[1,2]{Eysan Sans\thanks{sans20@itu.edu.tr}}
\author[1]{Melisa Akdemir\thanks{akdemirm20@itu.edu.tr}}
\author[1]{Ayse Tiryakioglu\thanks{aysetiryakioglu@itu.edu.tr}}
\author[1]{Ayse Peker-Dobie\thanks{pdobie@itu.edu.tr}}
\author[1]{Cihangir Ozemir\thanks{Corresponding author, ozemir@itu.edu.tr}}
\affil[1]{Department of Mathematics, Istanbul Technical University, Istanbul, Türkiye}
\affil[2]{Tobbund Logistics Investment Inc., Istanbul, Türkiye}

\begin{document}

\date{}

\maketitle

\begin{abstract}
Goodwin’s model is a cornerstone in the study of dynamical systems within macroeconomics, explaining the interaction between employment ratio and wage share in a closed economy. Analogous to predator-prey dynamics in mathematical economics, the Goodwin model, despite its simplicity, effectively captures the periodic behavior of state variables over specific time intervals. By relaxing the initial assumptions, the model can be adapted to account for more complex economic scenarios. In this article, we study a higher-dimensional extension of the Goodwin model that incorporates variable capacity utilization and capital coefficient alongside employment ratio and wage share. In particular instances, the wage share and employment rate equations decouple from the overall system. For these cases,  by incorporating a delay effect in the Phillips curve, we demonstrate that while the equilibrium of the generalized system remains stable within certain parameter domains in the absence of delay, the introduction of delay can induce a Hopf bifurcation, leading to periodic oscillations.  We analytically derive the critical delay parameter value that destabilizes the equilibrium point via a Hopf bifurcation.  

\vspace{0.5cm}
\end{abstract}

\section{Introduction}
The Goodwin model, originally proposed in \cite{goodwin1967growth}, predicts cyclic dynamics between the employment rate and wage share in a closed economy. In this article, we aim to study an extension  of this dynamic framework, starting with the following generalized form of the Goodwin model, as studied in \cite{glombowski1987generalizations} in the form
 \begin{subequations}\label{genmod}
\begin{eqnarray}
\dot{\beta} &=&  [\varphi(\beta) + f(\lambda) - u(\lambda) - n]\beta,\\
\dot{\lambda}&=&[\varphi(\beta)+ \phi_1(\beta)- u(\lambda) -\phi_2(\theta)]\lambda,  \\
\dot{\theta} &=& [\psi (\lambda) + f(\lambda) - z(\lambda,v)]\theta,\\
\dot{v} &=& [-f(\lambda)+z(\lambda,v)]v.
\end{eqnarray}
\end{subequations}
In this system, $\beta$ denotes the employment ratio, $\lambda$ represents the wage share, $\theta$ refers to the rate of capacity utilization, and $v$ is the actual capital coefficient. For brevity, we describe the functions on the right-hand sides of these equations in Section 2, where the analysis is presented. Initially, we examine the following delayed subsystem of this dynamical system
\begin{subequations}\label{main31int}
\begin{eqnarray}
\dot{\beta}(t) &=&  \Big[\beta_0+\gamma_2 \beta(t)-\delta_0 \lambda(t)  \Big]\beta(t), \label{delay31-1}\\
\dot{\lambda}(t)&=&\Big[ \lambda_0-\nu_2 \lambda(t)+\gamma_2 \beta(t)+\rho_1 \beta(t-\tau)  \Big]\lambda(t)  \label{delay32-4}
\end{eqnarray}
\end{subequations}
where $\beta_0$, $\gamma_2$, $\delta_0$, $\lambda_0$, $\nu_2$, and $\rho_1$ are constants.  We establish the conditions on the parameters under which the nonzero equilibrium point is stable or undergoes a Hopf bifurcation, subject to the delay  $\tau\geq 0$.  Furthermore, we conduct an analysis to determine the  direction of the Hopf bifurcation.  Subsequently, we perform a similar analysis for the following subsystem, which differs slightly from \eqref{main31int} by incorporating a new parameter and imposing an additional condition on the nonzero equilibrium point derived from the main system \eqref{genmod} as
\begin{subequations}\label{main32int}
\begin{eqnarray}
\dot{\beta}(t) &=&  \Big[\tilde\beta_0-\tilde \delta_0 \lambda(t)  \Big]\beta(t), \label{delay31-2}\\
\dot{\lambda}(t)&=&\Big[ \tilde \lambda_0-\mu_2\nu_2 \lambda(t)+\rho_1 \beta(t-\tau)  \Big]\lambda(t).  \label{delay32-1}
\end{eqnarray}
\end{subequations}
\subsection{Literature}
Since Goodwin introduced his simple yet crucial model on the dynamics between employment rate and wage share, inspired by the Lotka-Volterra predator-prey model, numerous extended models have been explored in subsequent studies. In addition to examining the stability of the original model, significant efforts have been directed towards enhancing its  realism. In this regard, some of these studies will be briefly explained.

Goodwin's original work \cite{goodwin1967growth} establishes the relationship between the employment ratio and the wage share in a closed economy through the following system of non-linear differential equations
\begin{subequations}
\begin{eqnarray}
\dot{v} &=& \Big[\frac{1}{\sigma}-(\alpha+\beta)-\frac{1}{\sigma }  u\Big] v, \\
\dot{u}&=&\Big[-(\alpha+\gamma) + \rho v \Big]u.
\end{eqnarray}
\end{subequations}
Extensions of Goodwin’s model, aiming to modify dynamics to keep state variable outputs within the unit square, are presented in \cite{desai2006clarification} and \cite{harvie2007dynamical}. Reference \cite{gaudenzi2013generalizations} follows a similar approach, applying the theory of planar Hamiltonian systems to examine Goodwin's class struggle model. One of the foundational assumptions in Goodwin’s model, namely that workers do not save, is relaxed in \cite{sordi2001growth}, where the share of capital held by workers is introduced as an additional variable. Considering variable capacity utilization and actual capital coefficient, \cite{glombowski1987generalizations} develops the four-dimensional system  \eqref{genmod} encompassing Goodwin's model as a special case. Building on Minsky’s Financial Instability Hypothesis \cite{minsky1982can}, Keen  integrates a banking sector into the Goodwin growth cycle in \cite{keen1995finance}, potentially disrupting the stable limit cycles in the long run and introducing instabilities, as analysed in \cite{grasselli2012analysis}. Extending this analysis, \cite{grasselli2015inflation} introduces inflation dynamics and speculative money flow into a four-dimensional system. In \cite{tebaldi2007chaotic}, modifications to Goodwin’s model include changes in technical progress assumptions and the introduction of a memory variable, leading to a three-dimensional system with a Hopf bifurcation observed at a critical memory effect value. Within the context of delayed dynamical systems, \cite{sportelli2022goodwin} considers finite time delays between investment orders and deliveries of finished capital goods, along with delayed reactions of real wages to unemployment levels. Reference \cite{desai1973growth} presents three extensions of Goodwin's model covering inflation, expected inflation, and excess capacity. Reference \cite{van1984implications} extends Goodwin's original model to incorporate differential savings decisions. In \cite{balducci1984generalization}, a dynamic game is formulated to explore cyclic behavior in Goodwin's model. \cite{ambrosi2015goodwin} proposes a new utility function for capitalists and workers, while \cite{harvie2000testing} empirically tests the Goodwin model across ten OECD countries. Comments on Harvie's work are discussed in \cite{grasselliharvie} and expanded upon in \cite{grasselli2018testing}, which introduce a capital accumulation rate between $(0,1]$  to improve Harvie's result. \cite{araujo2021testing} modifies the original Goodwin model by incorporating the rate of capacity utilization and validates it with US economic data. In \cite{yoshida2007dynamic}, the effects of government policy lag on macroeconomic variables are analysed using a three-dimensional system.
 
Finally, we revisit the model \eqref{genmod} as studied in \cite{glombowski1987generalizations}, where the authors assess the system's stability in four distinct scenarios. The first of these scenarios, labeled as \emph{ "a simplified Goodwin-like case"} \cite{glombowski1987generalizations}, represents the system \eqref{genmod} as follows
\begin{subequations}\label{genmod1}
\begin{eqnarray}
\dot{\beta} &=&  [f(\lambda) - \nu_1 - n]\beta,\\
\dot{\lambda}&=&[ \phi_1(\beta)- \nu_1]\lambda,  \\
\dot{v} &=& [-f(\lambda)+z(\lambda,v)]v,\\
\dot{\theta} &=& \frac{1}{\sigma v}.
\end{eqnarray}
\end{subequations}
Second, they examine the scenario of \emph{"variable speed of technical progress and changes in work intensity"} \cite{glombowski1987generalizations}, specifically focusing on the following subcase 
\begin{subequations}\label{genmod2}
\begin{eqnarray}
\dot{\beta} &=&  [\varphi(\beta) + f(\lambda) - \nu_1-\nu_2\lambda - n]\beta, \label{eqbeta1}\\
\dot{\lambda}&=&[\varphi(\beta)+ \phi_1(\beta)- \nu_1-\nu_2\lambda ]\lambda, \label{eqlambda1}\\
\dot{v} &=& [-f(\lambda)+z(\lambda,v)]v,\\
\theta &=& \frac{1}{\sigma v}.
\end{eqnarray}
\end{subequations}
Subsequently, they investigate another subcase under the title \emph{"non-neutral technical progress"} \cite{glombowski1987generalizations}, defined by 
\begin{subequations}\label{genmod3}
\begin{eqnarray}
\dot{\beta} &=&  [ f(\lambda) - u(\lambda) - n]\beta, \label{eqbeta2} \\
\dot{\lambda}&=&[\phi_1(\beta)-u(\lambda) ]\lambda, \label{eqlambda2}\\
\dot{v} &=& [-f(\lambda)+z(\lambda,v)]v, \label{eqv2}\\
\dot\theta&=& [\psi(\lambda)+f(\lambda)-z(\lambda,v)]\theta. \label{eqtheta2}
\end{eqnarray}
\end{subequations}
Finally, they perform the stability analysis of the complete system \eqref{genmod}. Systems \eqref{genmod2} and \eqref{genmod3} have been the primary models analysed in this work within the framework of delayed dynamical systems and Hopf bifurcation.
\subsection{Problem setup}
In the original Goodwin model, the Phillips curve illustrates the relationship between the inflation rate and (un)employment. Upon examining Phillips' seminal paper \cite{phillips1958relation}, particularly on page 297, it is evident that the curve derived from the data fitting process exhibits certain loops. Phillips notes that "\emph{A loop in this direction could result from a time lag in the adjustment of wage rates. If the rate of change of wage rates during each calendar year is related to unemployment lagged seven months, ... the loop disappears ... points lie closely along a smooth curve}" \cite{phillips1958relation}. This observation implies that the relationship depicted by the Phillips curve may  involve a time delay. This insight serves as our primary motivation for considering a generalized Goodwin model within the context of time delays in this article. 

In the derivation of system \eqref{genmod} presented in \cite{glombowski1987generalizations}, it becomes apparent that incorporating a time delay in the Phillips curve requires considering the function $\phi_1(\beta(t))$ with a delay, specifically in the form of $\phi_1(\beta(t-\tau))$. Thus, the analysis conducted in this article aims to address the following questions:
\begin{itemize}
\item Based on the results presented in \cite{glombowski1987generalizations},  system \eqref{genmod2} demonstrates stability at a nonzero equilibrium within a specific range of parameters. When extending this analysis to the delayed context described above, do parameter ranges exist such that the system (i) remains stable, (ii) undergoes a Hopf bifurcation, or (iii) becomes unstable?
\item What are the outcomes of applying the same analysis to system \eqref{genmod3}?
\end{itemize}
In Section 2.1, we explore the conditions necessary for a Hopf bifurcation in the context of a subsystem of \eqref{genmod2} by incorporating a time delay into the function $\phi_1$. The investigation focuses on two key criteria: (i) the presence of purely imaginary eigenvalues and (ii) the satisfaction of the transversality condition. Detailed calculations for the center manifold reduction are provided in the Appendix. We demonstrate that these conditions are satisfied for a specific parameter set. Section 2.2 performs  the same analysis for  a subsystem of \eqref{genmod3}, examining its Hopf bifurcation characteristics. 
As established in the literature, time delays can induce destabilization at equilibrium points in dynamical systems. We propose that the bifurcation and instability results obtained in this study may provide valuable insights into the potential mechanisms behind the instability of state variables or cyclic behaviors observed in certain systems, like those seen in \cite{harvie2000testing}. While these results may not offer quantitative predictions, they could still contribute qualitatively to understanding such phenomena.

\section{Analysis}
We will now detail the functions presented in \eqref{genmod}, as described in \cite{glombowski1987generalizations}
\begin{subequations}\label{funclist}
\begin{eqnarray}
\varphi(\beta)&=&-\gamma_1+\gamma_2 \beta, \label{fvarphi}\\
f(\lambda)&=&\delta [g(1-\lambda)-s_w],\label{ff}\\
u(\lambda)&=& \mu_1 + \mu_2 \nu_1+ \mu_2 \nu_2 \lambda,\label{fu}\\
\phi_1(\beta)&=&-\rho_0+\rho_1 \beta,\label{fa}\label{fphi1}\\
\phi_2(\theta) &=& b_2(1-a_3)\theta/(1-a_3b_3),\label{fphi2}\\
\psi(\lambda)&=& \nu_1 - \mu_1 - \mu_2 \nu_1 + \nu_2(1-\mu_2) \lambda,\label{fpsi}\\
z(\lambda,v)&=& c(1-\lambda)/v \label{fz}
\end{eqnarray}
\end{subequations}
where the constants are defined by
\begin{subequations}
\begin{eqnarray}
\rho_0&=&\dfrac{a_1 (1-b_3)-b_1 (1-a_3)}{1-a_3 b_3},\\
\rho_1&=&\dfrac{a_2 (1-b_3)}{1-a_3 b_3}.
\end{eqnarray}
\end{subequations}
$\varphi(\beta)$ is the work intensity depending on the employment rate, and the actual production rate is denoted by $f(\lambda)$. The technical capital productivity rate is given by $u(\lambda)$, while $\phi_1(\beta)$ and $\phi_2(\theta)$, depending on employment and capacity utilization respectively, represent functions for the real wage growth rate. $\psi(\lambda)$ refers to the technical capital coefficient rate, and $z(\lambda,v)$ represents the capital stock rate. The constraints for the constants, as specified in \cite{glombowski1987generalizations}, are outlined in Table 1.
\begin{table}[h]
\begin{center}
\caption{Restrictions on the parameters}\label{tab1}%
\begin{tabular}{@{}llll@{}}
\hline
$\mu_1 \ge 0 , \quad 1 \ge \mu_2 > 0$&  $b_1 \geq  0, \quad  b_2 \geq  0 , \quad  1 > b_3 \ge 0$\\
                    $\nu_1 \ge 0 , \quad  \nu_2 > 0 $& $1 > c > 0$ \\
                    $n \ge 0$ &  $1 > s_\pi > s_w > 0$ \\
                    $\gamma_1 \ge 0 , \quad  \gamma_2 \ge 0$ &  $\delta > 0$  \\
                    $a_1 > 0, \quad a_2 > 0 , \quad  1 \ge a_3 \ge 0$ & $1 > a_3 b_3$ \\
                    $g = c - (s_\pi - s_w)> 0$\\
\hline
\end{tabular}
\end{center}
\end{table}

\subsection{Variable speed of technical progress and changes in work intensity} 
According to the analysis in \cite{glombowski1987generalizations}, by utilizing variable speed of neutral technical progress and variable work intensity,  system \eqref{genmod} reduces to
\begin{subequations}\label{21main}
\begin{eqnarray}
\dot{\beta} &=&  [\varphi(\beta) + f(\lambda) - u(\lambda) - n]\beta, \label{eqbeta1}\\
\dot{\lambda}&=&[\varphi(\beta)+ \phi_1(\beta)-u(\lambda) ]\lambda, \label{eqlambda1}\\
\dot{v} &=& [-f(\lambda)+z(\lambda,v)]v,\\
\theta &=& \frac{1}{\sigma v}
\end{eqnarray}
\end{subequations}
where $\sigma $ is a constant.
In this case, the parameters are chosen to satisfy $\mu_1=0$, $\mu_2=1$, $\nu_1>0$, and $\nu_2>0$.  Since the first two equations, \eqref{eqbeta1} and \eqref{eqlambda1}, form a decoupled model, our analysis will focus on this subsystem.  The functions employed are defined as given in \eqref{fvarphi}, \eqref{ff}, and \eqref{fphi1}, with a specific instance of \eqref{fu} formulated by  $u(\lambda)=\nu_1+\nu_2 \lambda$. In the subsystem defined by \eqref{eqbeta1} and \eqref{eqlambda1},  $\phi_1(\beta)$ represents the relationship between the inflation and  employment rates. As previously discussed,  the literature suggests that this relation may involve a time delay. Thus, we incorporate this  delay into $\phi_1(\beta)$ and analyse the resulting delayed system
\begin{subequations}\label{main1}
\begin{eqnarray}
\dot{\beta}(t) &=&  \Big[\varphi\big(\beta(t)\big) + f\big(\lambda(t)\big) - u\big(\lambda(t)\big) - n\Big]\beta(t), \label{eqbeta1delay}\\
\dot{\lambda}(t)&=&\Big[\varphi\big(\beta(t)\big)+ \phi_1\big(\beta(t-\tau)\big)-u\big(\lambda(t))\Big]\lambda(t).   \label{eqlambda1delay}
\end{eqnarray}
\end{subequations}
If the functions $\varphi$, $f$, $u$, and $\phi_1$ are substituted, the system takes the following form
\begin{subequations}\label{main31}
\begin{eqnarray}
\dot{\beta}(t) &=&  \Big[\beta_0+\gamma_2 \beta(t)-\delta_0 \lambda(t)  \Big]\beta(t), \label{delay31-3}\\
\dot{\lambda}(t)&=&\Big[ \lambda_0-\nu_2 \lambda(t)+\gamma_2 \beta(t)+\rho_1 \beta(t-\tau)  \Big]\lambda(t)  \label{delay32-2}
\end{eqnarray}
\end{subequations}
where
\begin{align}
\beta_0&=(g-s_w)\delta -\gamma_1-\nu_1-n,\\
\lambda_0&=-\gamma_1-\rho_0-\nu_1,\\
\delta_0&=\nu_2+g\delta .
\end{align}

System \eqref{main31} represents a delayed logistic-type Lotka-Volterra system. While there is a substantial body of literature on delayed versions of Lotka-Volterra systems, we have not encountered any work that specifically addresses \eqref{main31}. Therefore, to the best of our knowledge, this system has not been previously investigated in the literature within the context we considered here.

We will focus on the equilibrium, defined by  
\begin{align}
\beta_e&=\frac{g\delta  (\gamma_1+\nu_1+\nu_2+\rho_0)-\nu_2(n+s_w \delta-\rho_0)}{\rho_1 \nu_2+g\delta (\rho_1+\gamma_2)},\\
\lambda_e&=\frac{\gamma_2 (g\delta+\rho_0-n-s_w \delta)-\rho_1 (n+\gamma_1+s_w \delta+\nu_1-g\delta)}{\rho_1 \nu_2+g\delta (\rho_1+\gamma_2)}
\end{align}
for  $ \beta_e, \lambda_e \in (0,1)$. The corresponding Jacobian matrices at the equilibrium $(\beta_e,\lambda_e)$ are given by
\begin{equation} 
J_0=\begin{bmatrix}
  \gamma_2 \beta_e   &    -\delta_0 \beta_e\\
\gamma_2 \lambda_e &-\nu_2 \lambda_e\\
\end{bmatrix}, \qquad \qquad 
J_\tau=\begin{bmatrix}
  0   &    0\\
\rho_1  \lambda_e &0\\
\end{bmatrix}.
\end{equation}
The characteristic equation for an eigenvalue of the linearization of  system \eqref{main31} is given by
\begin{equation}
|J_0+e^{-x\tau} J_{\tau}-xI|
=\begin{vmatrix}
  \gamma_2 \beta_e -x  &    -\delta_0 \beta_e\\
\gamma_2 \lambda_e +e^{-x\tau}\rho_1 \lambda_e &-\nu_2 \lambda_e-x\\
\end{vmatrix}=0,
\end{equation}
which yields
\begin{equation}\label{eigen}
P(x)=x^2+p_0 x+r_0+q_0 e^{-x\tau}=0
\end{equation}
where $p_0$, $r_0$ and $q_0$ are as follows
\begin{subequations}
\begin{eqnarray}
p_0&=&\nu_2 \lambda_e-\gamma_2\beta_e,\\
r_0&=&g\delta\gamma_2\beta_e\lambda_e,\\
q_0&=&\delta_0\rho_1\beta_e\lambda_e.
\end{eqnarray}
\end{subequations}
When the time delay 
$\tau$ is set to zero, the characteristic equation \eqref{eigen} simplifies to 
\begin{equation}\label{eigenzero}
P_0(x)=x^2+p_0 x+r_0+q_0 =0.
\end{equation}
\begin{lemma} Given that $r_0+q_0>0$, if
\begin{equation}\label{RH}
p_0>0,
\end{equation}
then all roots of the characteristic equation \eqref{eigenzero} have negative real parts. Consequently, the equilibrium point is stable for $\tau=0$.
\end{lemma}

Let $x=i\omega$ be a root of the characteristic equation \eqref{eigen}, where $\omega>0$. Then, the following equations hold
\begin{subequations}
\begin{eqnarray}
p_0\omega &=&q_0 \sin (\omega\tau),\\
\omega^2-r_0&=&q_0 \cos(\omega \tau).
\end{eqnarray}
\end{subequations}
By squaring and summing both sides of these equations, we obtain 
\begin{equation}
\omega^4+(p_0^2-2r_0)\omega^2+r_0^2-q_0^2=0.
\end{equation}
Letting $\omega^2=z$ where $z>0$, it simplifies to 
\begin{equation}\label{eqz}
h(z)=z^2+(p_0^2-2r_0)z+r_0^2-q_0^2=0.
\end{equation}
To have a purely imaginary eigenvalue $x=i\omega$, the equation must have at least one positive solution  $\tilde z$. Defining the discriminant as 
\begin{equation}
\Delta=(p_0^2-2r_0)^2-4(r_0^2-q_0^2), 
\end{equation}
we derive the following results.
\begin{itemize}
\item[(H1)] $\Delta <0 $: No real $\tilde z$.
\item[(H2)] $\Delta=0$ and $2r_0-p_0^2\leq 0$: No positive $\tilde z$.
\item[(H3)] $\Delta=0$ and $2r_0-p_0^2>0$: One positive $\tilde z$. 
\item[(H4)] $r_0^2-q_0^2<0$: One positive $\tilde z$.
\item[(H5)] $\Delta>0$, $r_0^2-q_0^2>0$ and $2r_0-p_0^2<0$: No positive $\tilde z$.
\item[(H6)] $\Delta>0$, $r_0^2-q_0^2>0$ and $2r_0-p_0^2>0$: Two positive $\tilde z$.
\end{itemize}
Therefore, we have the following outcome. 

\begin{lemma}
Equation \eqref{eqz} has the following properties.

(A) In cases (H1), (H2), (H5), \eqref{eqz} has no root  $\tilde z>0$. 

(B) In cases (H3), (H4) and (H6), \eqref{eqz} has at least one root $\tilde z >0$.
\end{lemma}

Suppose there exists a positive solution to \eqref{eqz}. Without loss of generality, we can assume that there are two positive roots $\tilde z_1$ and $\tilde z_2$. Consequently, there are two corresponding values of $\omega>0$ given by
\begin{equation}
\omega_{1}=\sqrt{\tilde z_1}, \qquad \omega_{2}=\sqrt{\tilde z_2}.
\end{equation}
Critical values of $\tau$ at which \eqref{eigen} has a purely imaginary root 
 $i\omega$ are determined by   
\begin{equation}\label{eqtau}
\tau_{k}^j=\frac{1}{\omega_k} \Big[\cos^{-1} \left(\frac{\omega_k^2-r_0}{q_0}\right)+2j\pi\Big]
\end{equation}
for  $k=1,2$ and $j=0,1,2,...$. Let us define
\begin{equation}
\tau_0=\tau_{k_0}^0=\mathrm{min}(\tau_1^0,\tau_2^0), \qquad \omega_0=\omega_{k_0}, \qquad z_0=\omega_0^2.
\end{equation}
Now we note the following Lemma from \cite{ruan2003zeros}.

\begin{lemma}\label{ruan}
	Consider the exponential polynomial equation
	\begin{align}
		P(x,e^{-x \tau_{1}},\ldots, e^{-x \tau_{m}})&= \, x^{n}+p_{1}^{(0)} x^{n-1}+\ldots+p_{n-1}^{(0)} x+p_{n}^{(0)}  +[p_{1}^{(1)} x^{n-1}+\ldots +p_{n-1}^{(1)} x +p_{n}^{(1)}]e^{-x \tau_{1}}   \nonumber   \\
		&+\ldots   \\
		&+[p_{1}^{(m)} x^{n-1}+\ldots+p_{n-1}^{(m)} x +p_{n}^{(m)}]e^{-x \tau_{m}}=0,   \nonumber
	\end{align}
	where $\tau_{i} \geq 0$ $(i=1,2,...,m)$ and $p_{j}^{(i)}$ $(i=0,1,2,...,m; j=1,2,...,n)$ are constants. As the parameters $(\tau_{1},\tau_{2},...,\tau_{m})$ vary, the sum of the order of the zeros of $P(x, e^{-x \tau_{1}},...,e^{-x \tau_{m}})$ on the open right half-plane can change only if a zero appears on or crosses the imaginary axis.
\end{lemma}
Based on this Lemma, the following conclusions can be drawn.
\begin{corollary}
	Suppose that the condition \eqref{RH} holds. 
	\begin{itemize}
		\item[(i)]If any of the conditions detailed in Lemma 2(A) holds, the roots of the characteristic equation \eqref{eigen} will possess negative real parts for all $\tau\geq 0$. 
		\item[(ii)] If any of the conditions detailed in Lemma 2(B) holds, the roots of the characteristic equation \eqref{eigen} will possess negative real parts for all $\tau\in [0,\tau_0)$. 
	\end{itemize}
\end{corollary}
\begin{lemma}
	Assume that one of the conditions stated in Lemma 2(B) holds, and that $h'(z_{0})\neq0$. In this case,  $i\omega_{0}$
  is a simple (i.e., non-multiple) purely imaginary root of the characteristic equation \eqref{eigen} when  $\tau=\tau_{0}$. Additionally, the  transversality condition 
	\begin{equation}\label{transvers}
		\frac{d(\mathrm{Re}(x(\tau)))}{d\tau}\Big\vert_{\tau=\tau_{0}}\neq0
	\end{equation}
	holds, and the sign of $d(\mathrm{Re}(x(\tau)))/d\tau\vert_{\tau=\tau_{0}}$ is consistent with the sign of $h'(z_{0})$.
\end{lemma}
\noindent
To demonstrate that $i\omega_0$ is a simple purely imaginary root of the characteristic equation \eqref{eigen} when $\tau = \tau_0$, provided that $h'(z_0) \neq 0$, we first assume, for contradiction, that $i\omega_0$ is a multiple root of \eqref{eigen}. Consequently, the condition $\frac{dP(x)}{dx}\big\vert_{x=i\omega_0} = 0$ must hold. Given that $x = i\omega_0$ is an eigenvalue when $\tau = \tau_0$, we analyse the equations
\begin{align}
P(i\omega_0)=(i\omega_0)^2+p_0 i\omega_0+r_0+q_0 e^{-i\omega_0\tau_0}=0,\\
P'(i\omega_0)=2i\omega_0+p_0-\tau_0 q_0e^{-i\omega_0\tau_0}=0.
\end{align}
By eliminating the exponential term from these equations, we derive 
\begin{equation}
p_0-\tau_0\omega_0^2+r_0\tau_0+i\omega_0(2+p_0\tau_0)=0.
\end{equation}
For this equality to hold, both the real and imaginary parts must separately equal zero. This results in
 $2\omega_0^2+p_0^2-2r_0=0$, which implies that $h'(\omega_0^2)=h'(z_0)=0$. This result contradicts our assumption that $h'(z_0) \neq 0$. Therefore, $i\omega_0$ must be a simple root of the equation. Consequently, from this part of the lemma, we utilize the result
\begin{equation}\label{nonzero}
\dfrac{dP}{dx}\Big\vert_{x=i\omega_0}\neq 0.
\end{equation}
To establish the validity of the transversality condition \eqref{transvers}, we proceed by differentiating the characteristic  equation \eqref{eigen} with respect to $\tau$. This differentiation yields
\begin{equation}
\left(2x+p_0-\tau q_0 e^{-x\tau}\right)\frac{dx}{d\tau}=xq_0 e^{-x\tau},
\end{equation}
which can also be expressed by
\begin{equation}
\frac{dP}{dx}\frac{dx}{d\tau}=xq_0 e^{-x\tau}.
\end{equation}
Using the result from \eqref{nonzero}, we find 
\begin{equation}
\frac{dx}{d\tau}\Big\vert_{\tau=\tau_0, \, x=i\omega_0}=\frac{i\omega_0 q_0 e^{-i\omega_0 \tau_0}}{\frac{dP}{dx}\Big\vert_{x=i\omega_0}}.
\end{equation}
Since the denominator is nonzero near $x = i\omega_0$ as established by \eqref{nonzero}, we obtain 
\begin{equation}
\frac{dx}{d\tau}\Big\vert_{\tau=\tau_0, \,x=i\omega_0}=\frac{i\omega_0 q_0 }{(2i\omega_0+p_0)e^{i\omega_0\tau_0}-\tau_0 q_0 }.
\end{equation}
Evaluating this expression at $x=i\omega_0$, we find
\begin{align}
\frac{d(\mathrm{Re}(x(\tau)))}{d\tau}\Big\vert_{\tau=\tau_{0}}&=\mathrm{Re}\left(\frac{dx}{d\tau}\Bigr\rvert_{\tau=\tau_0, \,x=i\omega_0}\right)=\mathrm{Re}\left(\frac{i\omega_0 q_0 }{(2i\omega_0+p_0)e^{i\omega_0\tau_0}-\tau_0 q_0 }\right) \nonumber\\
&=\frac{\omega_0^2(2\omega_0^2+p_0^2-2r_0)}{D}=\frac{\omega_0^2}{D}h'(\omega_0^2),
\end{align}
where 
\begin{equation}
D=[p_0\cos (\omega_0\tau_0)-2\omega_0\sin(\omega_0 \tau_0)-q_0 \tau_0]^2+[p_0\sin(\omega_0\tau_0)+2\omega_0\cos(\omega_0\tau_0)]^2.
\end{equation}
It is evident that $D \neq 0$ due to \eqref{nonzero}. Thus, the transversality condition \eqref{transvers} holds.
\subsubsection{Direction of Hopf bifurcation}
In this subsection, we initiate our analysis by outlining the preliminary steps and presenting the results prior to stating the main theorem. The detailed calculations are presented in the Appendix. In the preceding analysis, we focused on  subsystem \eqref{main31} and examined the parameter conditions necessary to satisfy the first two requirements for a Hopf bifurcation: (i) the existence of purely imaginary eigenvalues and (ii) the transversality condition. Assuming these two conditions are met, we proceed to evaluate the first Lyapunov coefficient.

Initially, we introduce a change of variable in  system \eqref{main31} by setting $t = \tau \tilde{t}$. Given that $\dfrac{d}{dt} = \dfrac{1}{\tau} \dfrac{d}{d\tilde{t}}$, we subsequently derive the following  
\begin{subequations}
\begin{eqnarray}
\dfrac{d}{d \tilde t}\,\beta(\tau \tilde t) &=&  \tau \Big[\beta_0+\gamma_2 \beta(\tau \tilde t)-\delta_0 \lambda(\tau \tilde t)  \Big]\beta(\tau \tilde t), \\
\dfrac{d}{d\tilde t }\,\lambda(\tau \tilde t)&=&\tau\Big[ \lambda_0-\nu_2 \lambda(\tau \tilde t)+\gamma_2 \beta(\tau \tilde t)+\rho_1 \beta(\tau (\tilde t-1))  \Big]\lambda(\tau \tilde t).  
\end{eqnarray}
\end{subequations}
By substituting $\tilde{\beta}(\tilde{t}) = \beta(\tau \tilde{t})$ and $\tilde{\lambda}(\tilde{t}) = \lambda(\tau \tilde{t})$ into the system, we obtain
\begin{subequations}
\begin{eqnarray}
\dfrac{d}{d \tilde t}\,\tilde \beta( \tilde t) &=&  \tau \Big[\beta_0+\gamma_2 \tilde \beta(\tilde t)-\delta_0 \tilde \lambda( \tilde t)  \Big]\tilde \beta( \tilde t), \\
\dfrac{d}{d\tilde t }\,\tilde \lambda( \tilde t)&=&\tau\Big[ \lambda_0-\nu_2 \tilde \lambda(\tilde t)+\gamma_2 \tilde \beta( \tilde t)+\rho_1 \tilde \beta( \tilde t-1)  \Big]\tilde \lambda( \tilde t).  
\end{eqnarray}
\end{subequations}
This normalizes the delay $\tau$ to 1, introducing a factor of $\tau$ on the right-hand side. Omitting  tildes, we have the  following
\begin{subequations}\label{main34}
\begin{eqnarray}
\dot{\beta}(t) &=&  \tau \Big[\beta_0+\gamma_2 \beta(t)-\delta_0 \lambda(t)  \Big]\beta(t), \\
\dot{\lambda}(t)&=&\tau\Big[ \lambda_0-\nu_2 \lambda(t)+\gamma_2 \beta(t)+\rho_1 \beta( t-1)  \Big]\lambda(t). 
\end{eqnarray}
\end{subequations}
We translate  the equilibrium point $E(\beta_e, \lambda_e)$ to the origin by defining
\begin{equation}
u_1(t)=\beta(t)-\beta_e, \qquad u_2(t)=\lambda(t)-\lambda_e.
\end{equation}
This transformation yields 
\begin{subequations}\label{main34}
\begin{eqnarray}
\dot u_1(t)=\tau \Big[\gamma_2 \beta_e u_1(t)&-&\delta_0 \beta_e u_2(t) +\gamma_2 u_1^2(t)-\delta_0 u_1(t) u_2(t)\Big ],\\
\dot u_2(t)=\tau \Big[\gamma_2 \lambda_e u_1(t)&-&\nu_2\lambda_e  u_2(t) +\rho_1 \lambda_e u_1(t-1)  \nonumber \\
                   &+&\gamma_2 u_1(t) u_2(t) 
              -\nu_2 u_2^2 (t) +\rho_1 u_1(t-1) u_2(t)\Big].  
\end{eqnarray}
\end{subequations}
Based on this transformed system, we proceed to determine the direction of  Hopf bifurcation and the stability of  bifurcating periodic solutions using the normal form theory as described in \cite{hassard1981theory}. For further details on the center manifold analysis, we refer to \cite{balachandran2009delay}.  The most similar, though not identical, model to \eqref{main31} found in the literature is studied in \cite{song2005local}. Following the methodologies outlined in \cite{song2005local}, as well as those in \cite{ccelik2013hopf} and \cite{manjunath2014stability}, we compute the following values (for detailed calculations, see  Appendix A)
\begin{eqnarray}
    c_1(0)&=&\frac{i}{2\omega\tau_k}\Big ( g_{11}g_{20}-2|g_{11}|^2-\frac{|g_{02}|^2}{3}\Big ) + \frac{g_{21}}{2}, \label{c10}\\    
     \bar \mu_{2}&=& -\frac{\mathrm{Re}(c_1(0))}{\mathrm{Re}(\lambda'(\tau_k))},\\
    \beta_2&=& 2\mathrm{Re}\big(c_1(0)\big ).  \label{beta2}
\end{eqnarray}
These quantities characterize the behavior of bifurcating periodic solutions at the critical value $\tau_k$. Specifically, $\bar\mu_2$ indicates the direction of  Hopf bifurcation: if $\bar\mu_2 > 0$ ($\bar\mu_2 < 0$), the bifurcation is supercritical (subcritical) and  bifurcating periodic solutions  exist for $\tau > \tau_k$ ($\tau < \tau_k$). $\beta_2$ assesses the stability of  bifurcating periodic solutions: the solutions are stable (unstable) on the center manifold if $\beta_2 < 0$ ($\beta_2 > 0$) \cite{song2005local}.

\begin{theorem}
	Consider the positive equilibrium point $E = E(\beta_e, \lambda_e)$ of  system \eqref{main31} within the unit box, and assume that the condition \eqref{RH} in Lemma 1 holds. Let $\tau_0$ be defined as specified previously. 
	\begin{itemize}
		\item  If one of the conditions listed in Lemma 2(A) holds, then the equilibrium point $E$ is asymptotically stable for all $\tau \geq 0$.

  \item If one of the conditions listed in Lemma 2(B) holds and $h'(z_0) \neq 0$, then $E$ remains asymptotically stable for $\tau \in [0, \tau_0)$. For $\tau \in (\tau_0, \tau_1)$, where $\tau_1$ is the first value obtained from \eqref{eqtau} that exceeds $\tau_0$, the equilibrium point $E$ becomes unstable. Thus,  system \eqref{main31} undergoes a Hopf bifurcation at $E$ when $\tau = \tau_0$.
	\end{itemize}
\end{theorem}

For a numerical example, the parameters for the system are chosen as  $\nu_1 = 0.02$, $\nu_2 = 0.04$, $\gamma_1 = 0.01$, $\gamma_2 = 0.012$, $\delta = 4.2$, $c = 0.38$, $n = 0.01$, $s_\pi = 0.24$, $s_w = 0.04$, $a_1 = 0.9$, $a_2 = 1$, $a_3 = 0.99$, $b_1 = 1.9$, $b_2 = 0$, and $b_3 = 0.6$.  In this case, we calculate the equilibrium values as $\beta_e = 0.90$ and  $\lambda_e = 0.70$, and we find  $z_0=0.501343$, $\omega_0=0.708056$ and $h'(z_0)=0.991515$. The critical delay parameter is $\tau_0 = 0.0348488$. This gives  $c_1(0)=0.00132164 - 0.0136561 i$, which leads to $\bar{\mu}_2 < 0$, indicating a subcritical Hopf bifurcation. Additionally, since $\beta_2 > 0$, the bifurcating periodic solution is unstable.

These findings are illustrated and further supported by Figures 1 to 6, which depict the dynamics of the solutions with the specified parameters as the delay parameter varies around its critical value.

\begin{remark}
The findings of this subsection can be summarized as follows. The special case \eqref{21main} of  system \eqref{genmod}, as derived in \cite{glombowski1987generalizations}, includes the subsystems \eqref{eqbeta1} and \eqref{eqlambda1}. The analysis conducted in \cite{glombowski1987generalizations} provides stability criteria for these subsystems and establishes conditions under which the equilibrium point remains stable. Our results show with concrete reference to the literature that, when considering the dynamics with delays, as formulated in \eqref{eqbeta1delay} and \eqref{eqlambda1delay}, the stability of the equilibrium point is significantly affected by the delay parameter. While the equilibrium point may be asymptotically stable for a certain parameter set in the absence of delays, the presence of delays introduces the potential for periodic oscillations or instability. Thus, the equilibrium's stability, which might be predicted to be stable in the non-delayed case, can vary depending on the value of the delay parameter relative to its critical threshold, as outlined in Theorem 1. In other words, while the non-delayed analysis might predict that the equilibrium point is asymptotically stable for a given set of parameters, the introduction of delays into the internal dynamics can alter this outcome. For the same parameter values, delayed dynamics can lead to periodic oscillations or instability.
\end{remark}

\subsection{Non-neutral technical progress} 
Allowing non-neutral technical progress and neglecting changes in work intensity and influence of capacity utilization,  authors of \cite{glombowski1987generalizations}  derive the following system of equations from the general model  in \eqref{genmod}
\begin{subequations}
\begin{eqnarray}
\dot{\beta} &=&  [ f(\lambda) - u(\lambda) - n]\beta, \label{eqbeta2} \\
\dot{\lambda}&=&[\phi_1(\beta)-u(\lambda) ]\lambda, \label{eqlambda2}\\
\dot{v} &=& [-f(\lambda)+z(\lambda,v)]v, \label{eqv2}\\
\dot\theta&=& [\psi(\lambda)+f(\lambda)-z(\lambda,v)]\theta \label{eqtheta2}
\end{eqnarray}
\end{subequations}
where the functions on the right-hand sides are as defined in \eqref{funclist}, and the parameters are required to satisfy the  conditions
\begin{equation}
\mu_1 \geq 0, \quad  0 < \mu_2 < 1, \quad \nu_1 \geq 0, \quad \nu_2> 0
\end{equation}
and 
\begin{subequations}
\begin{eqnarray}
u(\lambda)  &=& \mu_1 +\mu_2\nu_1+ \mu_2\nu_2\lambda,\\
\psi (\lambda) &=&\nu_1 -\mu_1 -\mu_2\nu_1(1-\mu_2)\lambda. 
\end{eqnarray}
\end{subequations}
Following our previous approach, we now consider the subsystem governed by equations \eqref{eqbeta2} and \eqref{eqlambda2}. Notably, as highlighted in \cite{glombowski1987generalizations}, for the nonzero equilibrium point $(\tilde \beta_e, \tilde \lambda_e, \tilde v_e, \tilde \theta_e)$, when $\psi(\lambda)\not\equiv 0$ the condition $\psi(\tilde \lambda_e)=0$ emerges from \eqref{eqv2} and \eqref{eqtheta2}. In accordance with \cite{glombowski1987generalizations}, this specific value of $\tilde \lambda_e$ is denoted as $\tilde \lambda_e^*$ and is given by 
\begin{equation}
\tilde \lambda_e^* = \frac{\mu_1 - \nu_1(1-\mu_2)}{\nu_2(1-\mu_2)}.
\end{equation}
Furthermore, from \eqref{eqbeta2} and \eqref{eqlambda2}, one obtains 
\begin{align}
\tilde \lambda_e &= \frac{\delta(g - s_w) - n - \mu_1-\mu_2\nu_1}{\delta g +\mu_2\nu_2}, \label{lambdae2}\\
\tilde \beta_e&=\frac{a_1(1-b_3)-b_1(1-a_3) + (1-a_3b_3)(\mu_1+\mu_2\nu_1+\mu_2\nu_2\tilde \lambda_e)}{a_2(1-b_3)}. \label{betae2}
\end{align}
Clearly, $\tilde \lambda_e^*$ and $\tilde \lambda_e$ must be identical. Assuming that the parameters satisfy $\tilde \lambda_e^* = \tilde \lambda_e$, we analyse the dynamics between $\beta(t)$ and $\lambda(t)$ in \eqref{eqbeta2} and \eqref{eqlambda2} with a delay term, similar to the approach in the previous subsection. Hence, we consider the model 
\begin{subequations}\label{main2}
\begin{eqnarray}
\dot{\beta}(t) &=&  \Big[ f\big(\lambda(t)\big) - u\big(\lambda(t)\big) - n\Big]\beta(t), \\
\dot{\lambda}(t)&=&\Big[\phi_1\big(\beta(t-\tau)\big)-u\big(\lambda(t))\Big]\lambda(t).  \end{eqnarray}
\end{subequations}
\eqref{main2} is a special case of \eqref{main1} where $\varphi(\beta)=-\gamma_1+\gamma_2 \beta\equiv 0$. The distinction between \eqref{main2} and \eqref{main1} lies in the constants involved in the term $u(\lambda)$. Therefore, the bifurcation analysis previously conducted with $\gamma_1 = \gamma_2 = 0$ will be analogous, with minor adjustments to the coefficients. The results for this particular case are summarized as follows. The subsystem can be expressed as
\begin{subequations}\label{main32}
\begin{eqnarray}
\dot{\beta}(t) &=&  \Big[\tilde\beta_0-\tilde \delta_0 \lambda(t)  \Big]\beta(t), \label{delay31-4}\\
\dot{\lambda}(t)&=&\Big[ \tilde \lambda_0-\mu_2\nu_2 \lambda(t)+\rho_1 \beta(t-\tau)  \Big]\lambda(t)  \label{delay32-3}
\end{eqnarray}
\end{subequations}
 where
\begin{align}
\tilde \beta_0&=(g-s_w)\delta- \mu_1-\mu_2\nu_1-n,\\
\tilde \lambda_0&=-\rho_0-\mu_1-\mu_2\nu_1,\\
\tilde \delta_0&=g\delta+\mu_2\nu_2 .
\end{align}
 We can re-express the nonzero equilibrium points  given by \eqref{lambdae2} and \eqref{betae2}     as  
\begin{equation}
\tilde \lambda_e=\frac{\tilde \beta_0}{\tilde \delta_0}, \qquad 
\tilde\beta_e=\frac{\mu_2\nu_2\tilde\beta_0-\tilde\lambda_0\tilde\delta_0}{\rho_1\tilde \delta_0}.
\end{equation}
Here we assume $\tilde \lambda_e, \tilde \beta_e \in (0,1)$.
The related Jacobian matrices at the equilibrium $\tilde E(\tilde \beta_e,\tilde \lambda_e)$ are
\begin{equation} 
J_0=\begin{bmatrix}
  0   &    -\tilde\delta_0 \tilde\beta_e\\
0     &-\mu_2\nu_2 \tilde \lambda_e\\
\end{bmatrix}, \qquad \qquad 
J_\tau=\begin{bmatrix}
  0   &    0\\
\rho_1  \tilde \lambda_e &0\\
\end{bmatrix}.
\end{equation}
The characteristic equation for an eigenvalue of the linearization of  the system \eqref{main32}
\begin{equation}
|J_0+e^{-x\tau} J_{\tau}-xI|
=\begin{vmatrix}
   -x  &    -\tilde\delta_0 \tilde\beta_e\\
e^{-x\tau}\rho_1 \tilde\lambda_e &-\mu_2\nu_2 \tilde\lambda_e-x\\
\end{vmatrix}=0
\end{equation}
yields  
\begin{equation}\label{eigen2}
\tilde P(x)=x^2+\tilde  p_0 x+\tilde q_0 e^{-x\tau}=0
\end{equation}
where
\begin{subequations}
\begin{eqnarray}
\tilde p_0&=&\mu_2\nu_2 \tilde\lambda_e,\\
\tilde q_0&=&\tilde \delta_0\rho_1\tilde \beta_e\tilde \lambda_e.
\end{eqnarray}
\end{subequations}
When $\tau=0$, \eqref{eigen2} simplifies to the following form 
\begin{equation}\label{eigenzero2}
\tilde   P_0(x)=x^2+\tilde  p_0 x+\tilde  q_0 =0.
\end{equation}
\begin{lemma} Since $\tilde p_0 >0$ and $\tilde q_0 >0$, 
all  roots of the characteristic equation \eqref{eigenzero} have negative real parts, and thus the equilibrium point is stable at $\tau=0$.
\end{lemma}

Let $x=i\omega$ be a root of \eqref{eigen2} with $\omega>0$. This implies
\begin{subequations}
\begin{eqnarray}
\tilde p_0\omega =&=&\tilde q_0 \sin (\omega\tau),\\
\omega^2&=&\tilde q_0 \cos(\omega \tau).
\end{eqnarray}
\end{subequations}
By squaring both equations and adding them, we derive 
\begin{equation}
\omega^4+\tilde p_0^2\omega^2-\tilde q_0^2=0.
\end{equation}
Substituting $\omega^2=z>0$, we obtain 
\begin{equation}\label{eqz2}
\tilde   h(z)=z^2+\tilde  p_0^2z-\tilde  q_0^2=0.
\end{equation}
To ensure the existence of a purely imaginary eigenvalue $x=i\omega$, equation \eqref{eqz2} must possess at least one positive solution $\hat z$.  Given that $\tilde h(0)=-\tilde q_0^2 <0$, it follows that \eqref{eqz2} indeed has at least one positive solution $\hat z$. 
In the general case, let us assume the existence of two positive roots, $\hat z_1$ and $\hat z_2$. Consequently, there are two corresponding positive values of $ \omega>0$, denoted by 
\begin{equation}
\tilde\omega_{1}=\sqrt{\hat z_1}, \qquad \tilde\omega_{2}=\sqrt{\hat z_2}.
\end{equation}
To determine the critical values of $ \tau$ such that \eqref{eigen2} has a purely imaginary root $i\tilde \omega$, we compute  
\begin{equation}\label{eqtau2}
\tilde\tau_{k}^j=\frac{1}{\tilde\omega_k} \Big[\cos^{-1} \left(\frac{\tilde\omega_k^2}{\tilde q_0}\right)+2j\pi\Big]
\end{equation}
for  $k=1,2$ and $j=0,1,2,...$. Let us define
\begin{equation}
\tilde\tau_0=\tilde\tau_{k_0}^0=\mathrm{min}(\tilde\tau_1^0,\tilde\tau_2^0), \qquad \tilde\omega_0=\tilde\omega_{k_0}, \qquad \hat z_0=\tilde\omega_0^2.
\end{equation}
Based on this, we can draw the following conclusion.

\begin{corollary}
For all $\tau \in [0, \tilde \tau_0)$, all roots of \eqref{eigen2} have negative real parts.
\end{corollary}

\begin{lemma}
	 The eigenvalue $i\tilde\omega_{0}$ is a simple (i.e., non-multiple) purely imaginary root of  equation \eqref{eigen2} when $\tau=\tilde \tau_{0}$. Furthermore, the following condition holds  
	\begin{equation}\label{transvers2}
		\frac{d(\mathrm{Re}(x(\tau)))}{d\tau}\Big\vert_{\tau=\tilde \tau_{0}}> 0,
	\end{equation}
	which confirms that the transversality condition is satisfied.
\end{lemma}
\noindent

The proof of Lemma 6  mirrors that of Lemma 4 with $\tilde r_0 = 0$. Therefore, we omit it here to avoid redundancy.

\begin{theorem}
	Let $\tilde E=\tilde E(\tilde \beta_e,\tilde \lambda_e)$ be a positive equilibrium point of system \eqref{main32},  and let $\tilde\tau_0$ be defined as previously. The equilibrium point  $\tilde E$ is asymptotically stable for $\tau \in [0, \tilde\tau_0 )$. It becomes unstable when $\tau \in (\tilde\tau_0, \tilde\tau_1 )$,  where $\tilde\tau_1$ is the first value of $\tau$ determined from  \eqref{eqtau2} that exceeds $\tilde\tau_0$. System \eqref{main32} undergoes a Hopf bifurcation at  $\tilde E$ when $\tau = \tilde\tau_0$.
 \end{theorem}
 
 The validity of Theorem 2 is further supported by Figures 7 through 12, produced with the parameter combination 
$ \mu_1 =0.0186145$, $\mu_2 = 0.5$, $\nu_1 = 0.015$ , $\nu_2 = 0.03$, $\gamma_1 = 0$, $\gamma_2 = 0$, $\delta = 4$, $c = 0.4$, $n = 0.01$, $s_\pi = 0.24$, $s_w = 0.04$, $a_1 = 0.9$, $a_2 = 1$, $a_3 = 1$, $b_1 = 1.9$, $b_2 = 0$, and $b_3 = 0.6.$ For this parameter set, the equilibrium values are $\tilde \beta_e=0.937$ and $\tilde \lambda_e^*=\tilde \lambda_e=0.741$, and the critical delay is $\tilde \tau_0=0.0196383$.  

For this scenario, we have not performed a directional analysis, as it closely follows the methodology outlined in the previous subsection. Additionally, the conclusions in this subsection can be justified by the arguments presented in Remark 1.

\section{Conclusion}

The Goodwin model, analogous to predator-prey dynamics in mathematical economics, despite its simplicity, effectively explains the periodic behavior of state variables observed within certain time intervals. This model’s assumptions can be relaxed to incorporate more complex economic scenarios, leading to various modifications in the literature. In this study, we explored a generalized, higher-dimensional Goodwin model that includes variable capacity utilization and capital coefficient in addition to the employment ratio and wage share.

By incorporating a delay effect in the Phillips curve, we demonstrated that while the equilibrium point of the generalized system remains stable within a specific parameter domain in the non-delayed case, the delayed model can experience a Hopf bifurcation, resulting in periodic oscillations. In particular instances, the equations for wage share and employment rate decouple from the full system. For these cases, we analytically determined the critical delay parameter value that causes the stable equilibrium point to become unstable via a Hopf bifurcation.

Through this analysis, we provide a fresh perspective on the original Goodwin model and its generalizations in the context of delayed dynamics. One noteworthy observation is that employment-wage share cycles constructed using real data for a specific country's economy are valid for a particular time interval before transitioning to a different equilibrium point, as can be seen in \cite{harvie2000testing}. This shift may be attributed to two primary factors: changes in system parameters or, as our findings may suggest, the delayed realization of the system's dynamics. The delay in the system’s response can lead to the destabilization of the initial equilibrium, causing the phase curve to shift and evolve into a new cycle within the phase space. Consequently, these results underscore the importance of incorporating delayed dynamics into the analysis of the Goodwin model and its generalizations, as they significantly impact the understanding of equilibrium transitions and economic cycle behavior.

The analysis of delay dynamical systems is well-documented in the literature; however, the Goodwin model and its variants have not been thoroughly examined within this framework. This research gap served as the primary motivation for our study. Our findings indicate that incorporating delays can substantially alter the stability and dynamic behavior of economic models, providing deeper insights into economic cycles and stability.

Additionally, we believe that the dynamics of the Goodwin model and its higher-dimensional extensions under delay deserve further exploration. Future research should investigate other modifications and extensions of the Goodwin model within the framework of delayed differential equations. This could potentially uncover new behaviors and stability criteria that are not apparent in non-delayed models, thereby enriching our understanding of economic dynamics.

In conclusion, our study bridges a critical gap in the literature by extending the Goodwin model to include delay effects and analysing its implications. The insights gained from this work underscore the importance of considering delays in economic models, as they can lead to complex dynamics. We hope this study will inspire further research into the delayed dynamics of economic models, paving the way for a deeper understanding of economic fluctuations and stability.

\appendix

\section{Appendix: Direction of Hopf bifurcation}

This Appendix provides the detailed technical calculations from Section 2.1.1 that lead to the results \eqref{c10}-\eqref{beta2}. Despite the extensive study of two-state variable  delayed-type Lotka--Volterra  systems in the literature, the specific case we have examined, to the best of our knowledge, has not been analysed before. We resume our analysis with  system \eqref{main34} from Section 2.1.1, continuing from where we previously concluded
\begin{subequations}\label{main34Appendix}
\begin{eqnarray}
\dot u_1(t)=\tau \Big[\gamma_2 \beta_e u_1(t)&-&\delta_0 \beta_e u_2(t) +\gamma_2 u_1^2(t)-\delta_0 u_1(t) u_2(t)\Big ],\\
\dot u_2(t)=\tau \Big[\gamma_2 \lambda_e u_1(t)&-&\nu_2\lambda_e  u_2(t) +\rho_1 \lambda_e u_1(t-1)  \nonumber \\
                   &+&\gamma_2 u_1(t) u_2(t) 
              -\nu_2 u_2^2 (t) +\rho_1 u_1(t-1) u_2(t)\Big].   
\end{eqnarray}
\end{subequations}
To ensure clarity, we first present this system in the following matrix form
\begin{align} \label{main342}
\begin{bmatrix}
  \dot u_1(t)\\
\dot u_2(t) 
\end{bmatrix}&=\tau\begin{bmatrix}
  \gamma_2 \beta_e   &    -\delta_0 \beta_e\\
\gamma_2 \lambda_e &-\nu_2 \lambda_e\\
\end{bmatrix}\begin{bmatrix}
   u_1(t)\\
 u_2(t)
\end{bmatrix} 
+\tau\begin{bmatrix}
  0   &    0\\
\rho_1  \lambda_e &0
\end{bmatrix}\begin{bmatrix}
   u_{1}(t-1)\\
 u_{2}(t-1) 
\end{bmatrix}   \nonumber\\
&+
\tau\begin{bmatrix}
  \gamma_2 u_1^2(t)-\delta_0 u_1(t) u_2(t)   \\
\gamma_2 u_1(t) u_2(t)-\nu_2 u_2^2 (t) +\rho_1 u_1(t-1) u_2(t)
\end{bmatrix}.
\end{align}
Let us define 
\begin{equation}
u_t(\theta)=\begin{bmatrix}
   u_{1t}(\theta)\\
 u_{2t}(\theta)
\end{bmatrix}:=u(t+\theta)=\begin{bmatrix}
   u_1(t+\theta)\\
 u_2(t+\theta)
\end{bmatrix}, \quad u_t: [-1,0]\rightarrow \mathbb{R}^2.
\end{equation}
Then, \eqref{main342} becomes
\begin{align} \label{A4}
\begin{bmatrix}
  \dot u_1(t)\\
\dot u_2(t) 
\end{bmatrix}&=\tau\begin{bmatrix}
  \gamma_2 \beta_e   &    -\delta_0 \beta_e\\
\gamma_2 \lambda_e &-\nu_2 \lambda_e   \nonumber \\
\end{bmatrix}\begin{bmatrix}
   u_{1t}(0)\\
 u_{2t}(0)
\end{bmatrix} 
+\tau\begin{bmatrix}
  0   &    0\\
\rho_1  \lambda_e &0
\end{bmatrix}\begin{bmatrix}
   u_{1t}(-1)\\
 u_{2t}(-1) 
\end{bmatrix}\\
&+
\tau\begin{bmatrix}
  \gamma_2 u_{1t}^2(0)-\delta_0 u_{1t}(0) u_{2t}(0)  \\
\gamma_2 u_{1t}(0) u_{2t}(0)-\nu_2 u_{2t}^2 (0) +\rho_1 u_{1t}(-1) u_{2t}(0)
\end{bmatrix}. 
\end{align}
For notational convenience, we set 
\begin{equation}
(a,b)=\begin{bmatrix}
  & a    &  \\
&  b     &
\end{bmatrix}_{2\times 1}, \qquad  
(a,b)^T=\begin{bmatrix}
  &a  &  b & 
\end{bmatrix}_{1\times 2}.
\end{equation}
For $\phi(\theta) =(\phi_1(\theta),\phi_2(\theta))  \in C([-1,0],\mathbb{R}^2)$, we define the operator $L_\mu$ by 
\begin{equation}  \label{hopf3}
\begin{split}
L_\mu\phi =\tau\begin{bmatrix}
  \gamma_2 \beta_e   &    -\delta_0 \beta_e\\
\gamma_2 \lambda_e &-\nu_2 \lambda_e\\
\end{bmatrix}\begin{bmatrix}
   \phi_1(0)\\
 \phi_2(0)
\end{bmatrix} 
+\tau\begin{bmatrix}
  0   &    0\\
\rho_1  \lambda_e &0
\end{bmatrix}\begin{bmatrix}
   \phi_1(-1)\\
 \phi_2(-1)
\end{bmatrix}
\end{split}
\end{equation}
and
\begin{equation}  \label{hopf4}
\begin{split}
f( \phi,\mu) =
\tau\begin{bmatrix}
  f_{11}  \\
f_{12}
\end{bmatrix} 
=
\tau\begin{bmatrix}
  \gamma_2 \phi_1^2(0)-\delta_0 \phi_1(0) \phi_2(0)  \\
\gamma_2 \phi_1(0) \phi_2(0)-\nu_2 \phi_2^2 (0) +\rho_1 \phi_1(-1) \phi_2(0)
\end{bmatrix} 
\end{split}
\end{equation}
with $\tau=\tau_k+\mu$.

Therefore, we have 
\begin{equation}\label{ulf}
\dot u(t)=L_\mu(u_t)+f(u_t,\mu).
\end{equation}
We aim to express the given system as a functional differential equation on $C([-1,0], \mathbb{R}^2)$ in the form
\begin{equation} \label{ulf2}
\dot{u}_t = A(\mu)u_t + R(\mu)u_t,
\end{equation}
which is more convenient since it involves a single unknown vector $u_t$, as opposed to \eqref{ulf}, that incorporates both $u(t)$ and $u_t$. By the Riesz representation theorem, there exists a matrix $\eta(\theta, \mu)$  with components that are functions of  bounded variation  for $\theta \in [-1, 0]$, such that
\begin{equation} \label{hopf5}
L_\mu \phi = \int_{-1}^{0} d\eta(\theta, \mu) \phi(\theta)
\end{equation}
for $\phi \in C[-1,0]$. The function $\eta(\theta, \mu)$ can be chosen as
\begin{equation} \label{hopf6}
\eta(\theta, \mu) = \tau \begin{bmatrix}
\gamma_2 \beta_e & -\delta_0 \beta_e \\
\gamma_2 \lambda_e & -\nu_2 \lambda_e
\end{bmatrix} \delta(\theta)
+ \tau \begin{bmatrix}
0 & 0 \\
\rho_1 \lambda_e & 0
\end{bmatrix} \delta(\theta + 1)
\end{equation}
where $\delta(\theta)$ is the Dirac distribution.
For $\phi(\theta) \in  C^1([-1,0], R^2)$,  define
\begin{equation} \label{hopf7}
A(\mu) \phi(\theta) =
\begin{cases} 
   \dfrac{d\phi(\theta)}{d\theta} ,   &\theta \in [-1,0),\\  \\
   \displaystyle\int_{-1}^{0}d\eta(s,\mu)\phi(s)=L_\mu(\phi),  &\theta=0,
\end{cases}
\end{equation}
and
\begin{equation} \label{hopf8}
R(\mu) \phi(\theta) =
\begin{cases} 
   0,   &\theta \in [-1,0),\\
   f( \phi,\mu),  &\theta=0.
\end{cases}
\end{equation}
Since $\dfrac{du_t}{d\theta}=\dfrac{du_t}{dt}$, \eqref{main342} is equivalent to \eqref{ulf2}. 
Indeed, for $\theta \in [-1,0)$, \eqref{ulf2} reduces to the trivial identity $\dfrac{du_t}{d\theta}=\dfrac{du_t}{dt}$, and it yields \eqref{main342} when $\theta=0$.

For $\psi  \in C^1  ([0,1],(R^2)^*)$,  define the adjoint $A^*(0)$ of $A(0)$ as 
\begin{equation} \label{adjoint}
A^*(0)\psi(s) =
\begin{cases} 
   -\dfrac{d\psi(s)}{ds} ,   &s\in (0,1],\\  \\ 
   \displaystyle\int_{-1}^{0}d\eta^T(t,0)\psi(-t),  &s=0
\end{cases}
\end{equation}
where $T$ denotes the  transpose, 
and a bilinear inner product  given by
\begin{equation} \label{hopf11}
   \langle \psi, \phi \rangle = \overline\psi(0)\cdot\phi(0) - \int_{\theta=-1}^{0}\int_{\xi=0}^{\theta}\overline{\psi}^T(\xi-\theta)d\eta(\theta)\phi(\xi)d\xi
\end{equation}
where $\eta(\theta) = \eta(\theta,0)$ and a bar denotes complex conjugate. Consequently, $A(0)$ and $A^*(0)$ are adjoint operators. Moreover, $\pm i \omega \tau_k$ are eigenvalues of $A(0)$ and hence are also eigenvalues of $A^*$. Let $q(\theta)$ be the eigenvector of $A(0)$ corresponding to $i \omega \tau_k$, and let $q^*(s)$ be the eigenvector of $A^*(0)$ corresponding to $-i \omega \tau_k$; that is,
\begin{equation}
A(0)q(\theta )=i\omega \tau_k q(\theta), \qquad A^*(0)q^*(s )=-i\omega \tau_k q^*(s).
\end{equation}
We need to find $q(\theta)$ that satisfies the following
\begin{equation}
\begin{cases}
			\dfrac{dq(\theta)}{d\theta}= i\omega \tau_k q(\theta)& \theta\in[-1,0),\\
            L_{\mu=0}\big(q(\theta)\big)=i\omega\tau_k q(0), & \theta=0,
\end{cases}
\end{equation}
and $q^*(s)$ that meets the conditions required by \eqref{adjoint}. It is then straightforward to demonstrate that
\begin{equation}
  q(\theta) = \big( 1, \alpha\big) e^{i\omega \tau_k\theta}, \qquad \alpha= \frac{\gamma_2\beta_e-i\omega}{\delta_0 \beta_e}
\end{equation}
and
\begin{equation}
  q^*(s) = B\big( \alpha^*,1\big) e^{i\omega \tau_k  s}, \qquad \alpha^*= \frac{i\omega-\nu_2 \lambda_e}{\delta_0 \beta_e}. 
\end{equation}
The evaluation of the inner product yields 
\begin{eqnarray}
\langle  q^*(s), q(\theta)\rangle &=&  \overline{q^*}(0)^T q(0)-\int_{\theta=-1}^0\int_{\xi=0}^\theta \overline{q^*}(\xi-\theta)^T \, d\eta(\theta,0) q(\xi) d\xi   \nonumber \\
&=& \overline{B} (\overline{\alpha}^*,1)^T(1,\alpha) \nonumber \\
&-&\int_{\theta=-1}^0\int_{\xi=0}^\theta  \overline{B} (\overline{\alpha}^*,1)^T e^{-i\omega\tau_k(\xi-\theta)}
d\eta(\theta,0)(1,\alpha)e^{i\omega\tau_k\xi}d\xi 
  \nonumber\\
&=& \overline{B} (\overline{\alpha}^*+\alpha)-\overline{B}\int_{\theta=-1}^0\theta e^{i\omega\tau_k \theta}(\overline{\alpha}^*,1)^T  d\eta(\theta,0) (1,\alpha)
\nonumber \\
\end{eqnarray}
\begin{eqnarray}
&=& \overline{B} (\overline{\alpha}^*+\alpha)- \overline{B}\int_{\theta=-1}^0\theta e^{i\omega\tau_k \theta}(\overline{\alpha}^*,1)^T\Bigg\{\tau_k\begin{bmatrix}
  \gamma_2 \beta_e   &    -\delta_0 \beta_e\\
\gamma_2 \lambda_e &-\nu_2 \lambda_e\\
\end{bmatrix}\delta(\theta) +\tau_k\begin{bmatrix}
  0   &    0\\
\rho_1  \lambda_e &0
\end{bmatrix}\delta(\theta+1)\Bigg\}(1,\alpha)d\theta 
\nonumber\\
&=& \overline{B} (\overline{\alpha}^*+\alpha)
+\overline{B} e^{-i\omega\tau_k }(\overline{\alpha}^*,1)^T \tau_k\begin{bmatrix}
  0   &    0\\
\rho_1  \lambda_e &0
\end{bmatrix} (1,\alpha)
\nonumber \\
&=&\overline{B} (\overline{\alpha}^*+\alpha) +\tau_k \overline{B} e^{-i\omega\tau_k} \rho_1\lambda_e   \nonumber \\
&=& \overline{B} \Big(\overline{\alpha}^*+\alpha+\tau_k  e^{-i\omega\tau_k} \rho_1\lambda_e\Big).
\end{eqnarray}
Now, we choose 
\begin{equation}
\overline{B}=\frac{1}{\overline{\alpha}^*+\alpha+\tau_k  e^{-i\omega\tau_k}\rho_1\lambda_e}
\end{equation}
so that 
\begin{equation}
     \langle q^*(s), q(\theta) \rangle = 1.
\end{equation}
We first compute the coordinates to describe the center manifold $C_0$ at $\mu =0$. Let $u_t$ be the solution of \eqref{ulf2} when $\mu =0$. We define
\begin{equation} \label{zteq}
z(t) = \langle q^*, u_t \rangle, \qquad  W(t,\theta)=u_t-zq-\bar z \bar q = u_t(\theta)- 2\mathrm{Re}\{z(t)q(\theta)\}.
\end{equation}
On the center manifold $C_0$, we have
\begin{equation}
    W(t,\theta)=W(z(t),\bar{z}(t),\theta)
\end{equation}
where
\begin{equation}
    W(z,\bar{z}, \theta)= W_{20}(\theta)\frac{z^2}{2}+ W_{11}(\theta)z\bar{z}+W_{02}(\theta)\frac{\bar{z}^2}{2}+W_{30}(\theta)\frac{z^3}{6}+...,
\end{equation}
z and $\bar{z}$ are local coordinates for the center manifold $C_0$ in the direction of $q^*$ and $\overline{q}^*$. We note that $W$ is real if $u_t$ is real, and therefore, we consider only real solutions. For a solution $u_t \in C_0$ of \eqref{ulf2}, given that $\mu =0$, we derive
\begin{eqnarray}
\dot{z} &=& \langle q^*, \dot u_t \rangle \nonumber \\
        &=& \langle q^*, A(0)u_t+R(0)u_t \rangle  \nonumber\\
        &=& \langle q^*, A(0)u_t \rangle+\langle q^*,R(0)u_t\rangle  \nonumber\\
         &=& \langle A^*(0)q^*, u_t \rangle+\langle q^*,f(u_t,0)\rangle 
 \nonumber\\  
         &=&\langle -i\omega\tau_k q^*,u_t\rangle+\langle q^*,f\big(W(z,\bar z,\theta) +2\mathrm{Re}\{z(t)q(\theta)\},0\big)   \rangle  \nonumber\\
         &=&\langle -i\omega\tau_k q^*,u_t\rangle+\langle q^*,f\big(W(z,\bar z,\theta) +2\mathrm{Re}\{z(t)q(\theta)\},0\big)   \rangle  \nonumber\\
         &=& i\omega\tau_k \langle q^*,u_t\rangle+ \overline{q}^*(0)^T f_0(z,\bar z)   \nonumber \\ 
         &=& i\omega\tau_k z(t)+g(z,\bar z)
\end{eqnarray}
where  $g(z,\bar z)= \overline{q}^*(0)^T f_0(z,\bar z) $, and 
\begin{equation} \label{hopf17}
    g(z,\bar{z}) =g_{20}\frac{z^2}{2}+g_{11}z\bar{z} +g_{02}\frac{\bar z^2}{2}+g_{21}\frac{z^2\bar z}{2}+... .
\end{equation}
Let us note that 
\begin{equation}
    u_t(u_{1t}(\theta),u_{2t}(\theta))=W(t,\theta)+zq(\theta)+\bar z \bar q (\theta),
\end{equation}
and $q(\theta)=(1,\alpha)e^{i\omega\tau_k\theta}$, and, explicitly,
\begin{eqnarray}
u_{1t}(0)&=&z+\bar z + W_{20}^{(1)}(0)\frac{z^2}{2}+W_{11}^{(1)}(0)z \bar z + W_{02}^{(1)}(0)\frac{\bar z^2}{2}+\ldots \nonumber \\
u_{2t}(0)&=& z\alpha +  \bar z \bar \alpha + W_{20}^{(2)}(0)\frac{z^2}{2}+W_{11}^{(2)}(0)z \bar z +W_{02}^{(2)}(0)\frac{\bar z^2}{2}+\ldots  \\
u_{1t}(-1)&=&ze^{-i\omega \tau_k}+\bar z e^{i\omega \tau_k} + W_{20}^{(1)}(-1)\frac{z^2}{2}+W_{11}^{(1)}(-1)z \bar z + W_{02}^{(1)}(-1)\frac{\bar z^2}{2}+\ldots \nonumber \\
u_{2t}(-1)&=&z\alpha e^{-i\omega \tau_k}+\bar z \bar \alpha e^{i\omega \tau_k} + W_{20}^{(2)}(-1)\frac{z^2}{2}+W_{11}^{(2)}(-1)z \bar z + W_{02}^{(2)}(-1)\frac{\bar z^2}{2}+\ldots \nonumber
\end{eqnarray}
\begin{eqnarray} \label{gzztilde}
    g(z,\bar z) &=& \overline{q}^*(0)^T f_0(z,\bar z)=\overline{B} \tau_k (\overline \alpha^*,1)^T 
    \begin{bmatrix}
  & f^0_{11}    & \\
&  f^0_{12}     &
\end{bmatrix}\nonumber \\
&=&\overline{B} \tau_k (\overline \alpha^*,1)^T\begin{bmatrix}
  \gamma_2 u_{1t}^2(0)-\delta_0 u_{1t}(0) u_{2t}(0)  \\
\gamma_2 u_{1t}(0) u_{2t}(0)-\nu_2 u_{2t}^2 (0) +\rho_1 u_{1t}(-1) u_{2t}(0)
\end{bmatrix}. 
\end{eqnarray}
Substituting \eqref{hopf17} into the left side of \eqref{gzztilde} and comparing the coefficients, we obtain
\begin{eqnarray}
g_{20}&=& 2  \overline B \tau_k (\alpha \gamma_2+\overline \alpha^* \gamma_2-\alpha \bar \alpha^*\delta_0-\alpha^2 \nu_2+\alpha \rho_1e^{-i\omega\tau_k }), \nonumber \\
g_{11}&=& \overline B \tau_k (\alpha \gamma_2+\overline \alpha \gamma_2+2 \overline\alpha^*\gamma_2-\alpha \overline\alpha^*\delta_0-\overline{\alpha}\; \overline{\alpha}^*\delta_0-2\alpha\overline\alpha\nu_2+e^{i\omega\tau_k}\alpha\rho_1+e^{-i\omega\tau_k}\overline\alpha\rho_1),\nonumber\\
g_{02}&=& 2  \overline B \tau_k (\overline\alpha\gamma_2+\overline\alpha^*\gamma_2-\overline\alpha \; \overline\alpha^*\delta_0-\overline\alpha^2\nu_2+e^{i\omega\tau_k}\overline\alpha\rho_1),\nonumber \\
g_{21}&=&\overline B \tau_k \Big(2\alpha\gamma_2 W^{(1)}_{11}(0)+4\overline\alpha^* \gamma_2 W^{(1)}_{11}(0)+\overline\alpha\gamma_2 W^{(1)}_{20}(0)+2\overline\alpha^* \gamma_2 W^{(1)}_{20}(0)\nonumber \\
         &+& 2\gamma_2 W^{(2)}_{11}(0)+\gamma_2W^{(2)}_{20}(0)-2\alpha \overline\alpha^* \delta_0  W^{(1)}_{11}(0)-\overline\alpha \; \overline\alpha^* \delta_0  W^{(1)}_{20}(0)-2\overline\alpha^* \delta_0  W^{(2)}_{11}(0)\nonumber \\
         &-&\overline\alpha^* \delta_0  W^{(2)}_{20}(0)
         -4\alpha \nu_2 W^{(2)}_{11}(0)-2\overline \alpha \nu_2 W^{(2)}_{20}(0) +2\alpha \rho_1 W^{(1)}_{11}(-1)\nonumber \\
         &+&\overline\alpha \rho_1 W^{(1)}_{20}(-1)+2\rho_1 e^{-i\omega\tau_k}W^{(2)}_{11}(0)+\rho_1 e^{i\omega\tau_k}W^{(2)}_{20}(0)\Big). 
\end{eqnarray}
Since  $g_{21}$  contains $W_{20}(\theta)$ and $W_{11}(\theta)$, we still need to compute them.
Using \eqref{ulf2} and \eqref{zteq}, we get
\begin{align}\label{hpf1}
    \dot{W}&= \dot{u}_t-\dot{z}q-\dot{\bar z}\bar q \nonumber\\
           &= A(0)u_t+R(0)u_t-(i\omega\tau_k z+g)q-(-i\omega \tau_k \bar z+\bar g)\bar q\nonumber  \\
           &=A(0)u_t+R(0)u_t-zA(0)q-\bar z A(0) \bar q-gq-\bar g \bar q \nonumber \\
           &=A(0)(u_t-zq-\bar z \bar q)+R(0)u_t-2\mathrm{Re}\{gq\} \nonumber \\
           &=A(0)W+R(0)u_t-2\mathrm{Re}\{gq\} \nonumber \\
           &=
    \begin{cases} 
    A(0)W-2\mathrm{Re}\{\bar q^*(0)\cdot 
 f_0q(\theta)\} ,   &\theta \in [-1,0),\\
   A(0)W-2\mathrm{Re}\{\bar q^*(0)\cdot  f_0q(0)\}+f_0,  &\theta=0,
\end{cases}\nonumber \\
&\overset{def}{=}A(0)W + H(z,\bar z,\theta)
\end{align}
where
\begin{equation} \label{h1}
   H(z,\bar z,\theta)=H_{20}(\theta)\frac{z^2}{2}+H_{11}(\theta)z\bar z+H_{02}\theta\frac{\bar z^2}{2}...
\end{equation}
Using $\dot W = W_z \dot z+W_{\bar z } \dot{\bar z}$ together with Eqs.   \eqref{hpf1} and \eqref{h1}, yields
\begin{eqnarray} \label{h3}
    (A(0)-2i\omega \tau_k)W_{20}(\theta)= - H_{20}(\theta),  \nonumber\\
    A(0)W_{11}(\theta)=-H_{11}(\theta),\\
    (A(0)+2i\omega \tau_k)W_{02}(\theta)= - H_{02}(\theta). \nonumber    
\end{eqnarray}
From \eqref{hpf1}, we know that for $\theta \in [-1, 0)$, 
\begin{equation}
     H(z,\bar z,\theta)=-\bar q^*(0)^T f_0q(\theta)-q^*(0)^T\bar f_0 \bar q(\theta)=-g(z,\bar z)q(\theta)-\bar{g}(z,\bar z)\bar q(\theta).
\end{equation}
Comparing the coefficients with \eqref{h1} gives 
\begin{equation} \label{h2}
\begin{split}
    H_{20}(\theta) = -g_{20} q(\theta)-\bar g_{02}\bar q(\theta), \\
    H_{11}(\theta) = -g_{11} q(\theta) -\bar g_{11}\bar q(\theta),
\end{split}
\end{equation}
and we get 
\begin{equation}
   \dot{W}_{20}(\theta) = 2i\omega \tau_k W_{20}(\theta)+g_{20}q(\theta) +\bar g_{02}\bar q(\theta).
\end{equation}
Note that $q(\theta)=q(0)e^{i\omega\tau_k\theta}$, and hence
\begin{equation}\label{h6}
    W_{20}(\theta) = \frac{ig_{20}}{\omega\tau_k}q(0)e^{i\omega\tau_k\theta}+ \frac{i \bar g_{02}}{3\omega\tau_k}\bar q(0)e^{-i\omega\tau_k\theta} + E_1e^{2i\omega\tau_k\theta}
\end{equation}
where $E_1=(E_{11},E_{12})$ is a constant vector. Similarly, from \eqref{h3} and \eqref{h2}, we derive
\begin{equation} 
    \dot{W}_{11}(\theta)= g_{11}q(\theta)+ \bar g_{11} \bar q(\theta),
\end{equation}
and 
\begin{equation} \label{h7}
     W_{11}(\theta) =- \frac{ig_{11}}{\omega\tau_k}q(0)e^{i\omega\tau_k\theta}+ \frac{i \bar g_{11}}{\omega\tau_k}\bar q(0)e^{-i\omega\tau_k\theta} + E_2
\end{equation}
where  $E_2=(E_{21},E_{22})$ is a constant vector.
In the following, we will identify appropriate $E_1$ and $E_2$ in \eqref{h6} and \eqref{h7}, respectively. It follows that
\begin{equation} \label{h20}
    \int_{-1}^{0} d\eta(\theta)W_{20}(\theta) = 2i\omega\tau_k W_{20}(0)-H_{20}(0),\\
\end{equation}
\begin{equation} \label{h25}    
     \int_{-1}^{0} d\eta(\theta)W_{11}(\theta) = - H_{11}(0)
\end{equation}
where $\eta(\theta) = \eta(\theta,0)$.
From \eqref{hpf1}, we have
\begin{equation} \label{h15}
    H_{20}(0) = -g_{20}q(0) - \bar g_{02} \bar q(0) +2\tau_k \begin{pmatrix}
        \gamma_2-\alpha \delta_0 \\
        \gamma_2 \alpha-\alpha^2\nu_2+\rho_1\alpha e^{-i\omega\tau_k}
    \end{pmatrix}
\end{equation}
and 
\begin{equation}\label{h16}
    H_{11}(0) = -g_{11}q(0) - \bar g_{11} \bar q(0) +\tau_k \begin{pmatrix}
        2\gamma_2-\alpha\delta_0-\overline\alpha\delta_0 \\
        2\gamma_2\mathrm{Re}\{\alpha\}-2\nu_2|\alpha|^2
        +2\rho_1 \mathrm{Re}\{\alpha  e^{i\omega\tau_k}\}
    \end{pmatrix}.
\end{equation}
Substituting \eqref{h6} and \eqref{h15} into \eqref{h20}, we obtain
\begin{equation}
    \begin{pmatrix}
        2i\omega\tau_k I-\displaystyle\int_{-1}^{0}e^{2i\omega\tau\theta}d \eta(\theta)
    \end{pmatrix} E_1 =2\tau_k \begin{pmatrix}
        \gamma_2-\alpha \delta_0 \\
        \gamma_2 \alpha-\alpha^2\nu_2+\rho_1\alpha e^{-i\omega\tau_k}
    \end{pmatrix},
\end{equation}
which simplifies to
\begin{equation}
    \begin{pmatrix}
        2i\omega-\gamma_2 \beta_e &\delta_0 \lambda_e\\
        -\gamma_2 \lambda_e-\rho_1 \lambda_e e^{-2i\omega\tau_k} &2i \omega+\nu_2 \lambda_e
    \end{pmatrix} E_1 =\begin{pmatrix}
         2\gamma_2-2\alpha \delta_0 \\
        2\gamma_2 \alpha-2\alpha^2\nu_2+2\rho_1\alpha e^{-i\omega\tau_k}
    \end{pmatrix}.
\end{equation}
Similarly, substituting \eqref{h7} and \eqref{h16} into \eqref{h25}, we get
\begin{equation}\label{hp1}
    \int_{-1}^{0}d\eta(\theta)E_2= -\tau_k \begin{pmatrix}
        2\gamma_2-\alpha\delta_0-\overline\alpha\delta_0 \\
        2\gamma_2\mathrm{Re}\{\alpha\}-2\nu_2|\alpha|^2
        +2\rho_1 \mathrm{Re}\{\alpha  e^{i\omega\tau}\}
    \end{pmatrix}.
\end{equation}
This implies
\begin{equation} \label{hp2}
    \begin{pmatrix}
        \gamma_2 \beta_e&-\delta_0 \lambda_e\\
        (\gamma_2+\rho_1)\lambda_e&-\nu_2 \lambda_e
    \end{pmatrix}E_2= - \begin{pmatrix}
        2\gamma_2-\alpha\delta_0-\overline\alpha\delta_0 \\
        2\gamma_2\mathrm{Re}\{\alpha\}-2\nu_2|\alpha|^2
        +2\rho_1 \mathrm{Re}\{\alpha  e^{i\omega\tau}\}
    \end{pmatrix}.
\end{equation}
It follows from  \eqref{h6}, \eqref{h7} \eqref{hp1}, and \eqref{hp2} that $g_{21}$ can be expressed. Thus, we can compute the following values
\begin{eqnarray}
    c_1(0)&=&\frac{i}{2\omega\tau_k}\Big ( g_{11}g_{20}-2|g_{11}|^2-\frac{|g_{02}|^2}{3}\Big ) + \frac{g_{21}}{2},\\    
    \bar \mu_2&=& -\frac{\mathrm{Re}(c_1(0))}{\mathrm{Re}(\lambda'(\tau_k))},\\
    \beta_2&=& 2\mathrm{Re}\big(c_1(0)\big ),
\end{eqnarray}
which determine the characteristics of bifurcating periodic solutions at the critical value $\tau_k$. For a set of parameters, numerical value of $c_1(0)$ and the implications drawn from these quantities have been presented in Subsection 2.1.1.

\vspace{1 cm}


\begin{figure}[h]
\centering
    \includegraphics[scale=0.8]{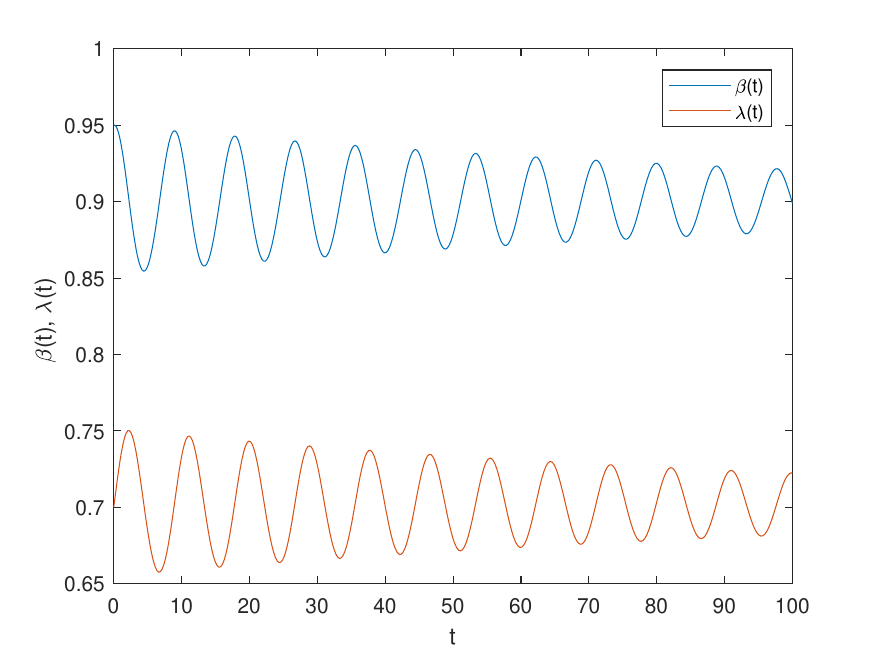}
    \caption{$\tau= 0 < \tau_0 = 0.0348488$, $t\in [0,100]$.}
\end{figure}
\begin{figure}[h]
\centering
    \includegraphics[scale=0.8]{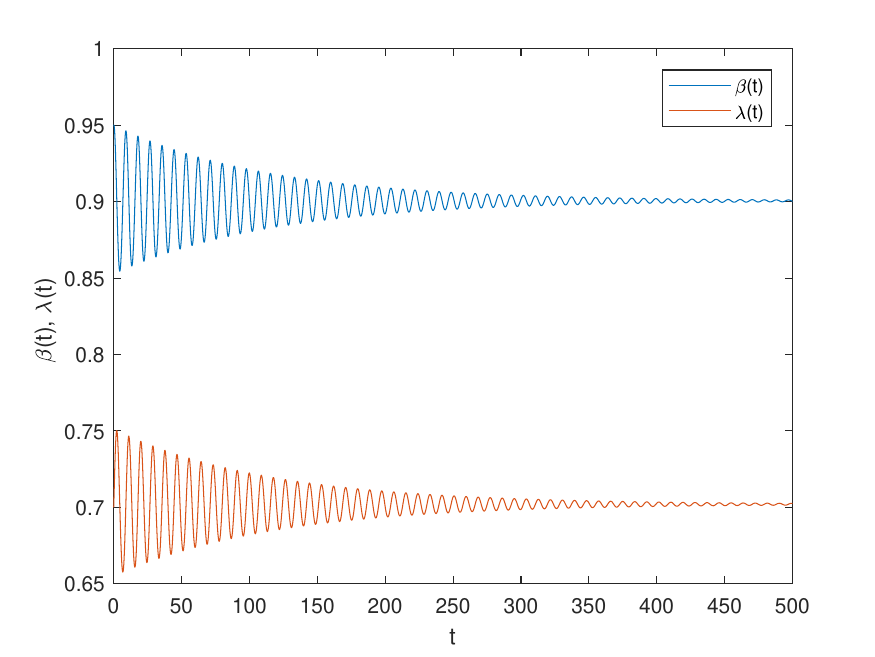}    
\caption{$\tau= 0 < \tau_0 = 0.0348488$, $t\in [0,500]$.}
\end{figure}
\begin{figure}[h]
\centering
    \includegraphics[scale=0.8]{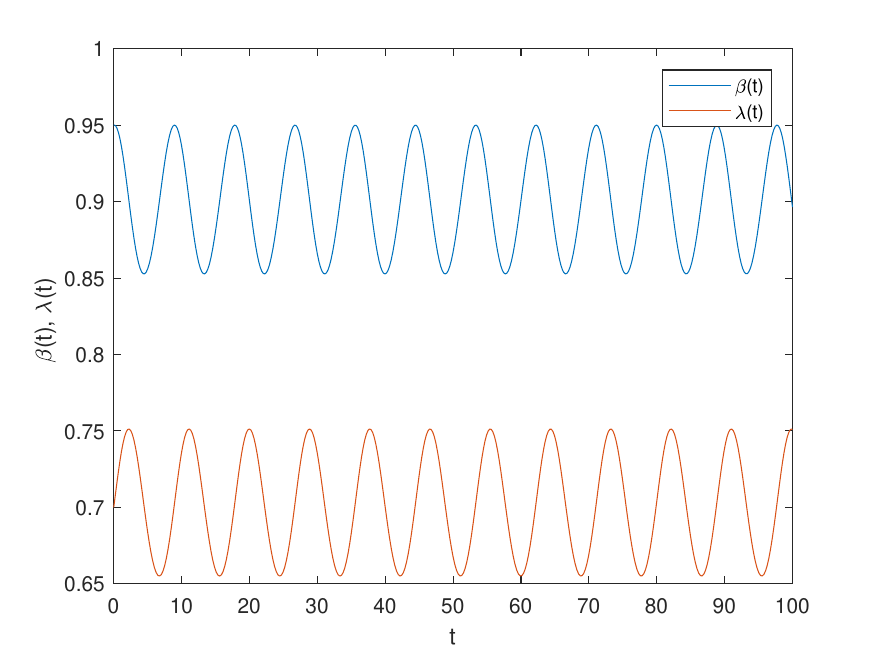}
\caption{$\tau = \tau_0   =0.0348488 $, $t\in [0,100]$.}
\end{figure}
\begin{figure}[h]
    \centering
    \includegraphics[scale=0.8]{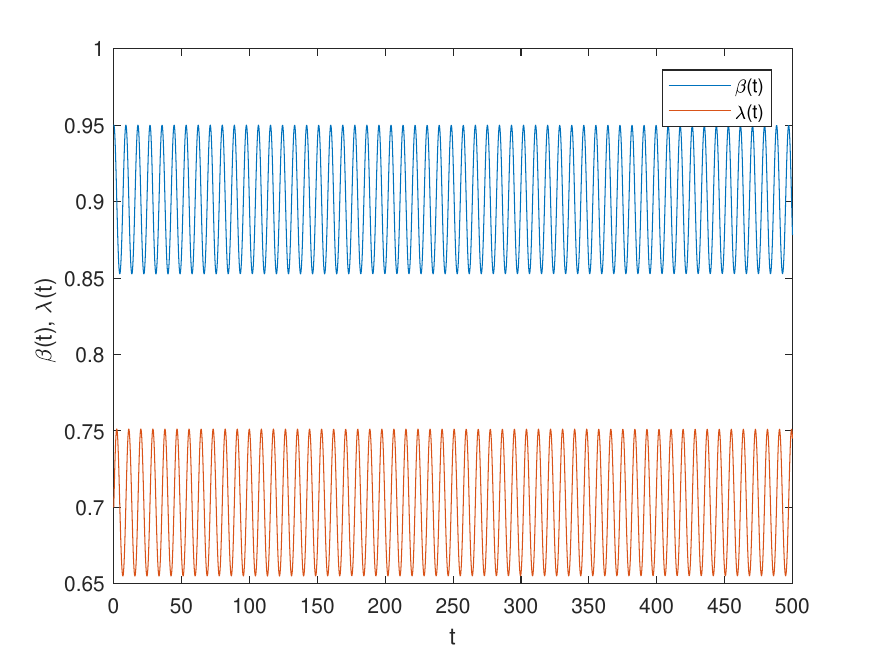}
\caption{$\tau = \tau_0   =0.0348488 $, $t\in [0,500]$.}
\end{figure}
\begin{figure}[h]
    \centering
    \includegraphics[scale=0.8]{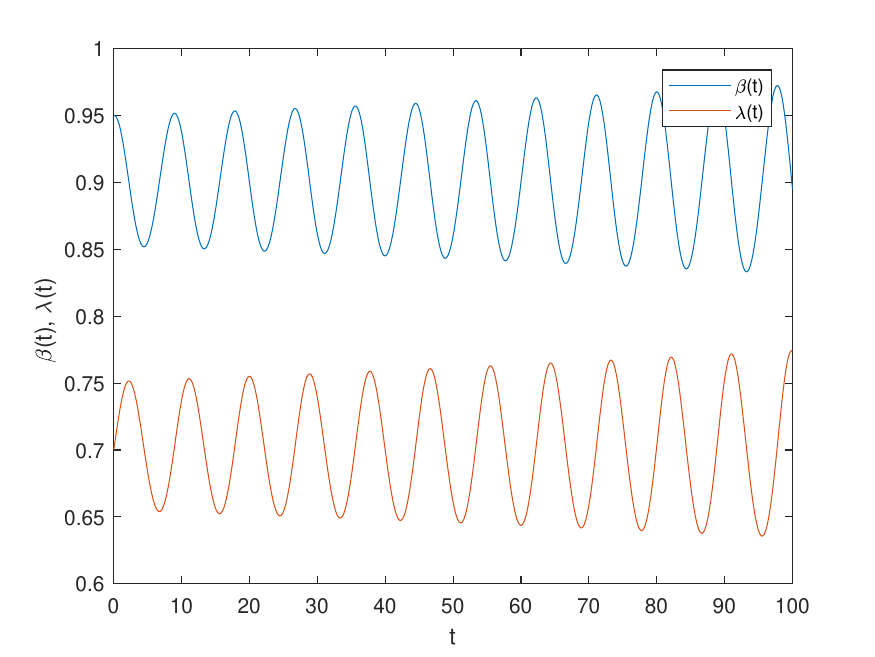}
    \caption{$\tau =0.05 > \tau_0  = 0.0348488$, $t\in [0,100]$.}\end{figure}
\begin{figure}[h]
    \centering
    \includegraphics[scale=0.8]{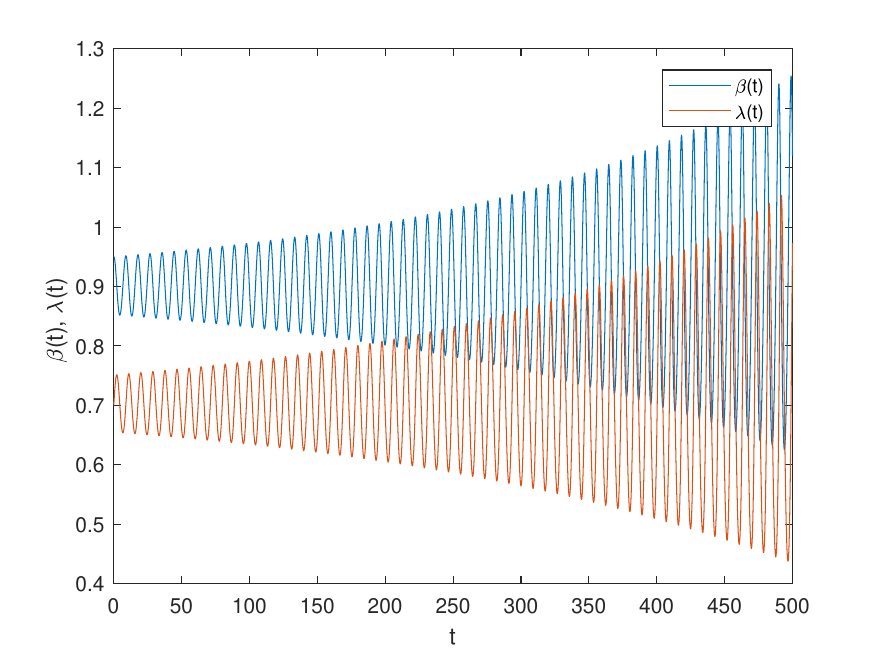}
\caption{$\tau =0.05 > \tau_0  = 0.0348488$, $t\in [0,100]$.}
\end{figure}

\begin{figure}[htbp]
\centering
\includegraphics[scale=0.8]{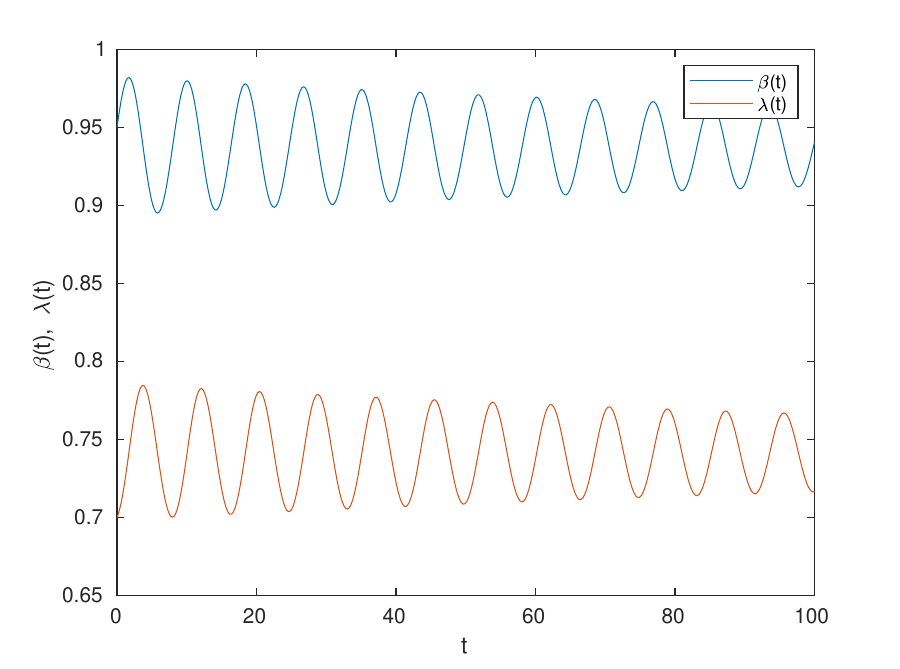}
\caption{$\tau  =0< \tilde\tau_0 = 0.0196383 $, $t\in [0,100]$.}\label{Figure3.4}
\end{figure}

\begin{figure}[htbp]
\centering
\includegraphics[scale=0.8]{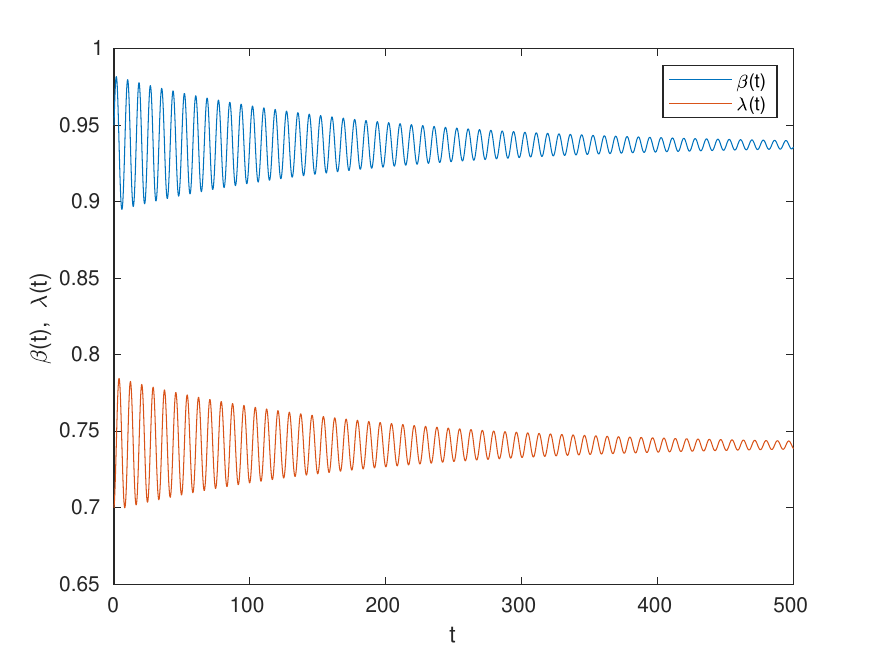}
\caption{$\tau  =0< \tilde\tau_0 = 0.0196383 $, $t\in [0,500]$.}\label{Figure3.5}
\end{figure}

\begin{figure}[htbp]
\centering
\includegraphics[scale=0.8]{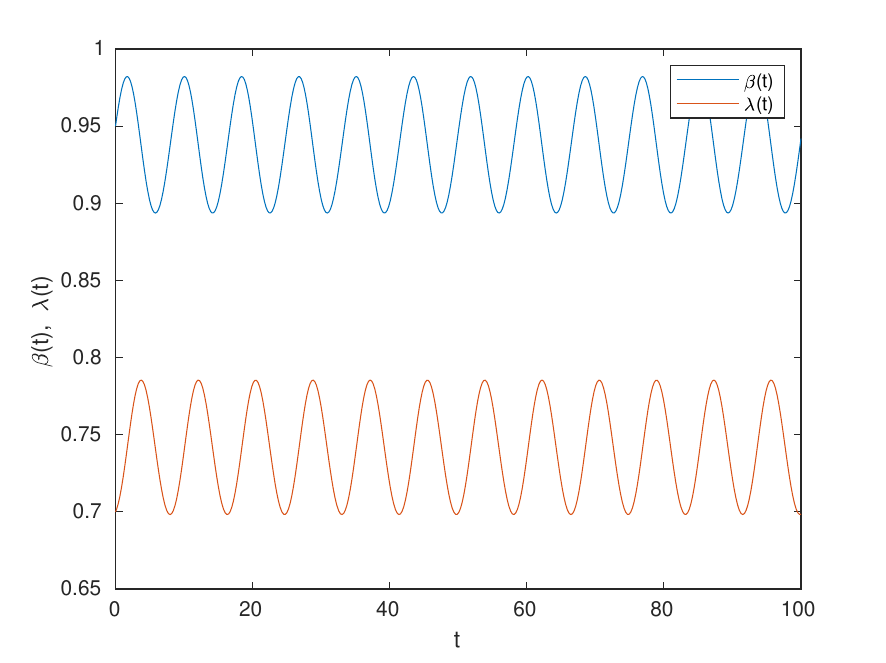}
\caption{$\tau = \tilde\tau_0 =0.0196383 $, $t\in [0,100]$.}\label{Figure3.8}
\end{figure}

\begin{figure}[htbp]
\centering
\includegraphics[scale=0.8]{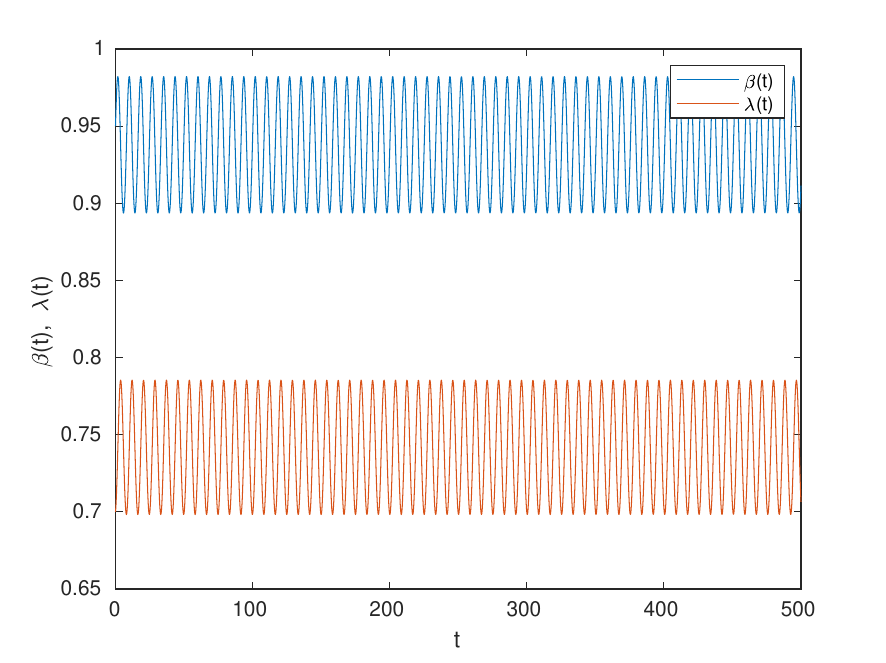}
\caption{$\tau= \tilde\tau_0 = 0.0196383$, $t\in [0,500]$.}\label{Figure3.9}
\end{figure}

\begin{figure}[htbp]
\centering
\includegraphics[scale=0.8]{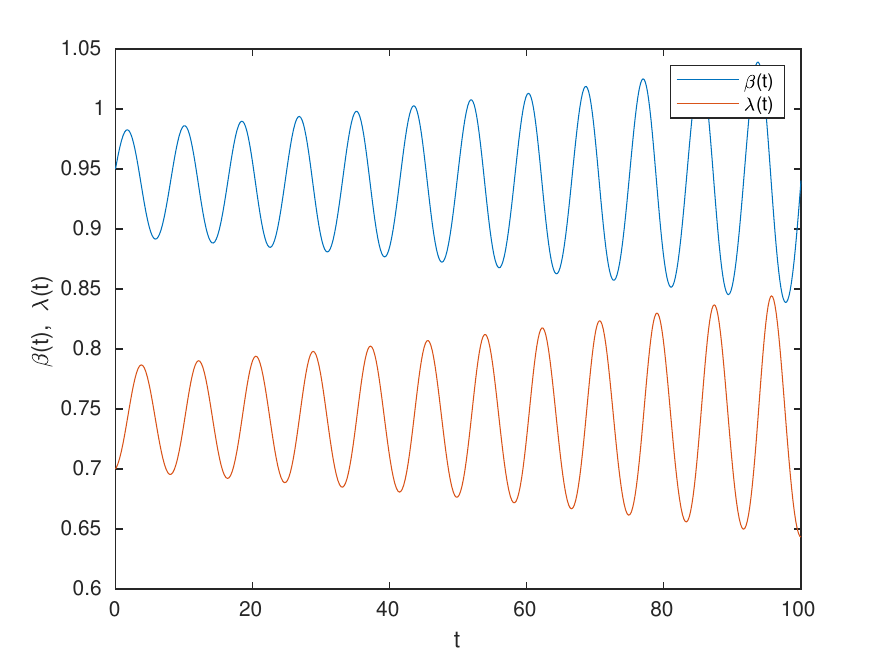}
\caption{$\tau = 0.05 > \tilde\tau_0  = 0.0196383 $, $t\in [0,100]$.}\label{Figure3.6}
\end{figure}

\begin{figure}[htbp]
\centering
\includegraphics[scale=0.8]{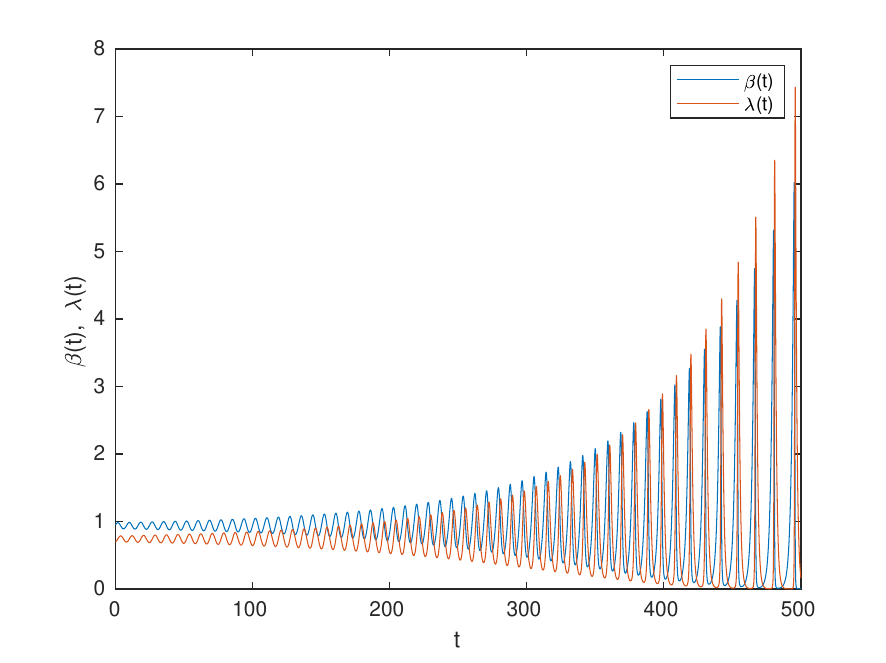}
\caption{$\tau  =0.05 > \tilde\tau_0 = 0.0196383 $, $t\in [0,500]$.}\label{Figure3.7}
\end{figure}

\end{document}